\documentclass[11pt]{article} \topmargin=-0.5cm
\usepackage{amssymb}

   \csname @addtoreset\endcsname{equation}{section}
    \csname @addtoreset\endcsname{figure}{section}
    \csname @addtoreset\endcsname{table}{section}

\newcounter{thanksnum}
\def\thanksnumber#1
{\setcounter{thanksnum}{\value{footnote}}\setcounter{footnote}{#1}%
                     \addtocounter{footnote}{-1}\footnotemark
                     \setcounter{footnote}{\value{thanksnum}}}
\def\newtheoremz#1{\@ifnextchar[{\@othmz{#1}}{\@nthmz{#1}}}

\def\@nthmz#1#2{%
\@ifnextchar[{\@xnthmz{#1}{#2}}{\@ynthmz{#1}{#2}}}

\def\@xnthmz#1#2[#3]{\expandafter\@ifdefinable\csname #1\endcsname
{\@definecounter{#1}\@addtoreset{#1}{#3}%
\expandafter\xdef\csname the#1\endcsname{\expandafter\noexpand
  \csname the#3\endcsname \@thmcountersepz \@thmcounterz{#1}}%
\global\@namedef{#1}{\@thmz{#1}{#2}}\global\@namedef{end#1}{\@endtheoremz}}}

\def\@ynthmz#1#2{\expandafter\@ifdefinable\csname #1\endcsname
{\@definecounter{#1}%
\expandafter\xdef\csname the#1\endcsname{\@thmcounterz{#1}}%
\global\@namedef{#1}{\@thm{#1}{#2}}\global\@namedef{end#1}{\@endtheoremz}}}

\def\@othmz#1[#2]#3{\expandafter\@ifdefinable\csname #1\endcsname
  {\global\@namedef{the#1}{\@nameuse{the#2}}%
\global\@namedef{#1}{\@thmz{#2}{#3}}%
\global\@namedef{end#1}{\@endtheoremz}}}

\def\@thmz#1#2{\refstepcounter
    {#1}\@ifnextchar[{\@ythmz{#1}{#2}}{\@xthmz{#1}{#2}}}

\def\@xthmz#1#2{\@begintheoremz{#2}{\csname the#1\endcsname}\ignorespaces}
\def\@ythmz#1#2[#3]{\@opargbegintheoremz{#2}{\csname
       the#1\endcsname}{#3}\ignorespaces}

\def\@thmcounterz#1{\noexpand\arabic{#1}}
\def\@thmcountersepz{.}
\def\@begintheoremz#1#2{ \trivlist \item[\hskip \labelsep{\bf #1\ #2}]}
\def\@opargbegintheoremz#1#2#3{ \trivlist
      \item[\hskip \labelsep{\bf #1\ #2\ (#3)}]}
\def\@endtheoremz{\endtrivlist}

\newtheorem{theorem}{Theorem}[section]
\newtheorem{lemma}{Lemma}[section]

\newtheorem{proposition}{Proposition}[section]

\newtheorem{condition}{Condition}[section]
\newtheorem{definition}{Definition}[section]
\newtheorem{remark}{Remark}[section]

\newtheorem{example}{Example}[section]

\def\e{\varepsilon}

\def\defi{\stackrel{{\scriptscriptstyle \Delta}}{=}}

\def\a{\alpha}
\def\d{\delta}
\def\o{\omega}
\def\O{\Omega}

\def\F{{\cal F}}
\def\w{\widehat}
\def\Ind{{\mathbb{I}}}

\def\esssup{\mathop{\rm ess\, sup}}
\def\essinf{\mathop{\rm ess\, inf}}

\def\dist{{\rm dist\,}}
\def\R{{\bf R}}
\def\E{{\bf E}}
\def\P{{\bf P}}

\def\Z{{\cal Z}}

\def\h{h}
\def\L{L}

\def\b{\beta}
\def\s{\delta}
\def\g{\gamma}

\def\W{{\cal W}^*}
\def\ww{\widetilde}
\def\X{{\cal X}}
\def\t{\theta}
\def\oo{\bar}
\def\s{\sigma}

\def\p{\partial}
\def\G{\Gamma}

\def\V{{\cal V}}
\def\A{{\cal A}}
\def\M{{\cal M}}

\def\L{{\cal L}}

\def\h{h}

\def\TT{{\cal T}}
\newcommand{\be}{\begin{equation}}
\newcommand{\ee}{\end{equation}}
\newcommand{\bd}{\begin{displaymath}}
\newcommand{\ed}{\end{displaymath}}
\newcommand{\ba}{\begin{array}{ll}}
\newcommand{\ea}{\end{array}}
\newcommand{\baa}{\begin{eqnarray}}
\newcommand{\eaa}{\end{eqnarray}}
\newcommand{\baaa}{\begin{eqnarray*}}
\newcommand{\eaaa}{\end{eqnarray*}}   \font\sm=cmr10

\def\ww{\tilde}
\def\W{{\cal W}}

\def\u{ u }

\def\Q{{\cal Q}}
\def\CC{{\cal C}}
\def\ttau{\tau_{{\scriptscriptstyle T}}}
\def\tttau{\tau_{{\scriptscriptstyle T},\d}}

\def\LU{{\scriptscriptstyle LU}}
\def\LS{{\scriptscriptstyle L}}
\def\US{{\scriptscriptstyle U}}
\def\dcoer{\varrho}
\date{\sm {Submitted: November 26, 2012. Revised: May 23, 2014}}
\title{\index{Degenerate backward SPDEs in domains with boundaries and with non-local conditions}
{Degenerate backward SPDEs in domains: non-local boundary
conditions and applications to finance}}
\author{
Nikolai Dokuchaev\\
 {\sm Department of Mathematics \& Statistics, Curtin
University,}\\ {\sm  GPO Box U1987, Perth, 6845 Western Australia} }
\begin{document}
\maketitle
\let\thefootnote\relax\footnote{This is a
pre-copy-editing, author-produced PDF. A part of this paper (entire Section 3 and the corresponding proofs from Section 6)
are included in an article "Degenerate backward SPDEs in bounded   domains and
applications to barrier options " accepted for
publication in {\em Discrete and Continuous Dynamical Systems --
Series A (DCDS-A)} following peer review. }

\begin{abstract}
Backward stochastic partial differential equations of parabolic type
in bounded domains are studied in the setting where the coercivity
condition is not necessary satisfied. Some generalized solutions
based on the representation theorem are suggested. For the backward
equation with a Cauchy condition at the terminal time,  some
regularity is derived from the regularity of the first exit times of
non-Markov characteristic processes. In addition, problems with
special non-local in time boundary conditions are considered. These
non-local conditions connect the terminal value of the solution with
a functional over the entire past solution. Uniqueness, solvability
and regularity results are obtained. Some applications to portfolio
selection problem are considered.
\end{abstract}
\section{Introduction}
Partial differential equations and  stochastic partial
differential equations (SPDEs) have fundamental significance for
natural sciences, and various boundary value problems for them
were widely studied.
 Usually,  well-posedness of a  boundary value depends on
the choice of the boundary value conditions.
\par
 Boundary value problems for SPDEs are well studied  in the existing literature
 for the case
 of  forward and backward parabolic Ito equations with the  Cauchy condition at
initial time or terminal time respectively (see, e.g., Al\'os et al
(1999), Bally {\it et al} (1994),  Da Prato and Tubaro (1996),
Gy\"ongy (1998), Krylov (1999), Maslowski (1995), Pardoux (1993),
 Rozovskii (1990), Walsh (1986), Zhou (1992), and Dokuchaev (1992), (2005), (2011), (2012)
and the bibliography there).  Many results have been also obtained
for the  pairs of   forward and backward
 equations with separate  Cauchy conditions at initial time and
 the terminal time respectively; see, e.g., Yong and Zhou
 (1999).
\par
Usually, SPDEs of parabolic types are considered under some
assumptions of coercivity such as Condition \ref{cond3.1.A} below
with $\dcoer >0$. Without this condition, an equation is regarded as
degenerate.  For the degenerate backward SPDEs in the whole space,
i.e., without boundaries, regularity results were obtained in
Rozovskii (1990),   Ma and Yong (1996, 1997), Hu {\em et al} (2002), Hamza and Klebaner (2005).
\index{
Rozovskii (1990),   Ma and Yong (1996, 1997), Hu {\em et al} (2002).}
In Rozovskii (1990),   Ma and Yong (1996, 1997), Hu {\em et al} (2002), second order parabolic type degenerate  SPDEs were considered.
In  Hamza and Klebaner (2005),  first order forward SPDEs were considered; these equations also can be classified
as degenerate.   For problems with boundaries, a different class of backward first order SPDEs was introduced in Bender and Douchaev (2014).

 The
methods applied in these works cannot be applied in the case of a
domain with a boundary  because of regularity issues that prevent
using of approximation of the differential operator  by a
non-degenerate one. It turns out that the theory of degenerate SPDEs
in domains is much harder than in the whole space and was not
addressed yet in the existing literature. Regularity is a difficult
issue for degenerate equations in the presence of a boundary.

We address this problem again.
The main contribution of this paper is an existence and regularity result for.
a parabolic type homogeneous
backward SPDEs that can be degenerate; the coercivity condition is
not necessary satisfied. We suggest a generalized solutions based on
the representation theorem for the backward equation with a Cauchy
condition at the terminal time.  Some regularity is obtained; the
proof is based on  the regularity of the first exit times of
non-Markov characteristic processes (Theorem \ref{ThRep}). The  result for the domains with a boundary is new. Our
proof is based on  the estimates from Dokuchaev (2004,2008a) of the $L_1$-distances between  the first exit times of
characteristic processes of underlying
backward SPDEs;  these estimates imply  estimates (\ref{tau}) and (\ref{tau2}) below that were crucial for the proof.
The estimates for non-Markov processes were established in  Dokuchaev (2008a) in a setting
that covered only the case of one dimensional
processes and the case of vector processes in domains with lacunas. Because of this,  we consider only these two types of the domains in non-Markov case (Condition \ref{condO}(ii)).
\par
The second contribution of this paper is extension on the case of degenerate SPDEs  of the earlier results
from Dokuchaerv (2013) for SPDEs with non-local conditions. In the literature,
there are many results for SPDEs with boundary conditions connecting the solution at different times, for instance, at the initial time and at the
terminal time. This category includes stationary type solutions for
forward SPDEs (see, e.g., Caraballo {\em et al } (2004),
 Chojnowska-Michalik (19987), Chojnowska-Michalik and Goldys
 (1995), Duan {\em et al} (2003), Mattingly (1999)
Mohammed  {\em et al} (2008), Sinai (1996), and the references
here).  Periodic solutions of SPDEs were also studied
 (Chojnowska-Michalik (1990), Feng and Zhao (2012), Kl\"unger  (2001)).
 As was mentioned in Feng and Zhao (2012), it is difficult to expect that, in general, a SPDE has a periodic
 in time solution $u(\cdot,t)|_{t\in[0,T]}$ in a usual sense of exact equality  $u(\cdot,t)=u(\cdot,T)$ that holds almost surely
 given that $u(\cdot,t)$ is adapted to some Brownian motion.
The periodicity of the solutions of stochastic equations was
usually considered
 in the sense of the distributions. In Feng and Zhao (2012),
 the periodicity was established  in a  stronger sense as a "random
 periodic solution (see Definition 1.1 from Feng and Zhao (2012)).
 Some periodic stochastic solutions were
obtained in Rodkina {\em et al} (2014) in some asymptotic sense for a
setting with time decaying  random noise. In Dokuchaev (2008b), the standard
boundary value Cauchy condition at the one fixed time was replaces by
a condition that mixes in one equation the initial value of the
solution and a functional of the
 entire solution for a forward SPDE.  In Dokuchaev (2013),
we considered non-local conditions that included almost surely periodicity for non-degenerate backward SPDEs; the prior estimates were obtained  in $L_2$-setting.  \par
The present paper addresses these and related problems with non-local boundary conditions again for degenerate backward SPDEs.  We
consider linear  Dirichlet  condition on the boundary of the state
domain. Instead of the Cauchy condition at the terminal time,  we consider conditions such as
 $\t^{-1}\int_0^\t u(\cdot,t)dt=u(\cdot,T)$ a.s., as well as  more general
conditions. Related results were obtained in Dokuchaev (2013)
for non-degenerate backward SPDEs  in $L_2$-setting. i..e,  with prior estimates
for the solutions based on $L_2$-norm. The novelty of the results of the present paper for SPDEs with non-local conditions, with respect to the related parer
of  Dokuchaev (2013),
  is that we allow equations to be degenerate; in addition, we consider
the prior estimates in $L_\infty$-setting. 
\par
 We present sufficient conditions for existence and
regularity of the solutions. As an example of
applications, a solution of portfolio selection problem is obtained
for continuous time market model with random coefficients  (Theorem \ref{ThF}).
\section{The problem setting and definitions}
 We
are given a standard  complete probability space $(\O,\F,\P)$ and a
right-continuous filtration $\F_t$ of complete $\s$-algebras of
events, $t\ge 0$. We assume that $\F_0$ is the $\P$-augmentation of
the set $\{\emptyset,\O\}$. We are given also a $N$-dimensional
Wiener process $w(t)$ with independent components;  it is a Wiener
process with respect to $\F_t$.
\par
Assume that we are given a bounded open domain $D\subset\R^n$  with a
$C^2$-smooth boundary $\p D$. Let $T>0$ be given, and let $Q\defi
D\times [0,T]$. \par
 We will study the following boundary value
problem in $Q$
\begin{eqnarray} 
\label{parab1} &&d_tu+(\A u+ \varphi)\,dt +\sum_{i=1}^N
B_i\chi_idt=\sum_{i=1}^N\chi_i(t)dw_i(t), \quad t\ge 0,
\\\label{parab10}
&& u(x,t,\o)\,|_{x\in \p D}=0
\\ &&u(\cdot, T)-\G u(\cdot)=\xi.
\label{parab2}
\end{eqnarray}
Here $u=u(x,t,\o)$, $\chi_i=\chi_i(x,t,\o)$,
$\varphi=\varphi(x,t,\o)$, $\xi=\xi(x,\o)$,
 $(x,t)\in Q$,   $\o\in\O$.\par
   In (\ref{parab2}), $\G$ is a linear operator that maps functions
defined on $Q\times \O$  to functions defines on $D\times \O$. For
instance, the case where  $\G u=u(\cdot,0)$ is not excluded; this
case corresponds to the periodic type boundary condition \baa
u(\cdot,T)-u(\cdot,0)=\xi.\label{period}\eaa
\par
  In (\ref{parab1}),  \baa \A
v=\sum_{i,j=1}^nb_{ij}(x,t,\o)\frac{\p^2 v}{\p x_i \p x_j}(x)
+\sum_{i=1}^n f_i(x,t,\o)\frac{\p v}{\p x_i
}(x)+\lambda(x,t,\o)v(x), \label{A}\eaa
 where $b_{ij}, f_i, x_i$ are the
components of $b$, $f$, and $x$ respectively,  and \be\label{B}
B_iv\defi\frac{dv}{dx}\,(x)\,\beta_i(x,t,\o),\quad i=1,\ldots ,N.
\ee
\par
We assume that the functions $b(x,t,\o):
\R^n\times[0,T]\times\O\to\R^{n\times n}$, $\b_j(x,t,\o):
\R^n\times[0,T]\times\O\to\R^n$, $f(x,t,\o):
\R^n\times[0,T]\times\O\to\R^n$, $\lambda(x,t,\o):
\R^n\times[0,T]\times\O\to\R$,  $\varphi (x,t,\o): \R^n\times
[0,T]\times\O\to\R$, and $\chi_i(x,t,\o): \R^n\times
[0,T]\times\O\to\R$,  are progressively measurable with respect to
$\F_t$ for all $x\in\R^n$, and the function $\xi(x,\o):
\R^n\times\O\to\R$ is $\F_0$-measurable for all $x\in\R^n$.
\subsection*{Spaces and classes of functions} 
We denote by $\|\cdot\|_{ X}$ the norm in a linear normed space
$X$, and
 $(\cdot, \cdot )_{ X}$ denote  the scalar product in  a Hilbert space $
X$.
\par
We introduce some spaces of real valued functions.
\par
 Let $G\subset \R^k$ be an open
domain, then ${W_q^m}(G)$ denote  the Sobolev  space of functions
that belong to $L_q(G)$ together with the distributional
derivatives up to the $m$th order, $q\ge 1$.
\par
 We denote  by $|\cdot|$ the Euclidean norm in $\R^k$, and $\bar G$ denote
the closure of a region $G\subset\R^k$.
\par Let $H^0\defi L_2(D)$,
and let $H^1\defi \stackrel{\scriptscriptstyle 0}{W_2^1}(D)$ be the
closure in the ${W}_2^1(D)$-norm of the set of all smooth functions
$u:D\to\R$ such that  $u|_{\p D}\equiv 0$. Let $H^2=W^2_2(D)\cap
H^1$ be the space equipped with the norm of $W_2^2(D)$. The spaces
$H^k$ and $W_2^k(D)$ are called  Sobolev spaces, they are Hilbert
spaces, and $H^k$ is a closed subspace of $W_2^k(D)$, $k=1,2$.
\par
 Let $H^{-1}$ be the dual space to $H^{1}$, with the
norm $\| \,\cdot\,\| _{H^{-1}}$ such that if $u \in H^{0}$ then
$\| u\|_{ H^{-1}}$ is the supremum of $(u,v)_{H^0}$ over all $v
\in H^1$ such that $\| v\|_{H^1} \le 1 $. $H^{-1}$ is a Hilbert
space.
\par We shall write $(u,v)_{H^0}$ for $u\in H^{-1}$
and $v\in H^1$, meaning the obvious extension of the bilinear form
from $u\in H^{0}$ and $v\in H^1$.
\par
We denote by $\oo\ell _{k}$ the Lebesgue measure in $\R^k$, and we
denote by $ \oo{{\cal B}}_{k}$ the $\sigma$-algebra of Lebesgue
sets in $\R^k$.
\par
We denote by $\oo{{\cal P}}$  the completion (with respect to the
measure $\oo\ell_1\times\P$) of the $\s$-algebra of subsets of
$[0,T]\times\O$, generated by functions that are progressively
measurable with respect to $\F_t$.
\par
 We  introduce the spaces
 \baaa
 &&X^{k}(s,t)\defi L^{2}\bigl([ s,t ]\times\Omega,
{\oo{\cal P }},\oo\ell_{1}\times\P;  H^{k}\bigr), \quad\\ &&Z^k_t
\defi L^2\bigl(\Omega,{\cal F}_t,\P; H^k\bigr),\\
&&\CC^{k}(s,t)\defi C\left([s,t]; Z^k_T\right), \qquad k=-1,0,1,2,
\\&& \X^k_c= L^{2}\bigl([ 0,T ]\times\O,\, \oo{{\cal P}
},\oo\ell_{1}\times\P;\; C^k(\oo D)\bigr),\quad \Z^k_c\defi
L_2(\O,\F_T,\P;C^k(D)),\quad k\ge 0. \eaaa
  The
spaces $X^k(s,t)$ and $Z_t^k$  are Hilbert spaces.
 \par
We introduce the spaces $$ Y^{k}(s,t)\defi
X^{k}(s,t)\!\cap \CC^{k-1}(s,t), \quad k=1,2, $$ with the norm $ \|
u\| _{Y^k(s,T)}
\defi \| u\| _{{X}^k(s,t)} +\| u\| _{\CC^{k-1}(s,t)}. $
For brevity, we shall use the notations
 $X^k\defi X^k(0,T)$, $\CC^k\defi \CC^k(0,T)$,
and  $Y^k\defi Y^k(0,T)$.

Let $\oo\V=L_{\infty}(\O,\F_T, C(\oo D))$, and let $\V$ be the set
of all $v\in \oo\V$  such that $v(x)|_{x\in\p D}=0$ a.s. We consider
$\V$ as a Banach space equipped with the norm of $\V$.

For a set ${\rm S}$ and a Banach space ${\rm X}$, we denote by ${\rm
B}({\rm S},{\rm X})$ the Banach space of bounded functions ${\rm x}:{\rm S}\to{\rm X}$
equipped with the norm $\|{\rm x}\|_{{\rm B}}=\sup_{{\rm s}\in {\rm S}}\|{\rm x}({\rm s})\|_{{\rm
X}}$.
\par
Let  $U={\rm B}([0,T];\V)\cap X^0\cap \CC^0$. We consider  $U$ as a
Banach space equipped with the norm ${\rm B}([0,T]; \V)$.
\par
Sometimes we shall omit $\o$.
\subsection*{Conditions on the domain and the coefficients}
 To proceed further, we assume that Conditions
\ref{cond3.1.A}-\ref{cond3.1.B} remain in force throughout this paper.
\begin{condition} \label{cond3.1.A} There exists a constant
$\dcoer \ge 0$ such that
\baaa
 \label{Main1} y^\top  b
(x,t,\o)\,y-\frac{1}{2}\sum_{i=1}^N |y^\top\b_i(x,t,\o)|^2 \ge
\dcoer |y|^2 \quad\forall\, y\in \R^n,\ (x,t)\in  D\times [0,T],\
\o\in\O. \eaaa
\end{condition}
\par
Condition \ref{cond3.1.A} with $\dcoer >0$ is usually called the
coercivity condition. If Condition \ref{cond3.1.A} is satisfied for
$\dcoer =0$ only, (\ref{parab1}) is usually referred as a  degenerate
equation. This important case  is included in this paper.
\begin{condition}\label{cond3.1.B}
The functions  $f(x,t,\o)$, $\lambda (x,t,\o)$, and $\b_i(x,t,\o)$
are bounded. These functions  are differentiable
in $x$ for a.e. $t,\o$, and the corresponding derivatives are
bounded. In addition, $b\in \X_c^3$, $\w f\in\X_c^2$,
$\lambda\in\X^1_c$, $ \b_i\in\X_c^3$, and $\b_i(x,t,\o)=0$ for
$x\in \p D$, $i=1,...,N$.
\end{condition}
\subsection*{The definition of solution}
\subsubsection*{Solution from $Y^1$}
\begin{proposition} 
\label{propL} Let $\zeta\in X^0$,
 let a sequence  $\{\zeta_k\}_{k=1}^{+\infty}\subset
L^{\infty}([0,T]\times\O, \ell_1\times\P;\,C(D))$ be such that all
$\zeta_k(\cdot,t,\o)$ are progressively measurable with respect to
$\F_t$, and  $\|\zeta-\zeta_k\|_{X^0}\to 0$  as $k\to +\infty$. Let $t\in [0,T]$ and
$j\in\{1,\ldots, N\}$ be given.
 Then the sequence of the
integrals $\int_0^t\zeta_k(x,s,\o)\,dw_j(s)$ converges in $Z_t^0$ as
$k\to\infty$, and its limit depends on $\zeta$, but does not depend
on $\{\zeta_k\}$.
\end{proposition}
\par
{\it Proof} follows from completeness of  $X^0$ and from the
equality
\begin{eqnarray*}
\E\int_0^t\|\zeta_{k}(\cdot,s,\o)-\zeta_m(\cdot,s,\o)\|_{H^0}^2\,ds
=\int_D\,dx\,\E\left(\int_0^t\big(\zeta_k(x,s,\o)-
\zeta_m(x,s,\o)\big)\,dw_j(s)\right)^2.
\end{eqnarray*}
\begin{definition} 
\rm Let $\zeta\in X^0$, $t\in [0,T]$, $j\in\{1,\ldots, N\}$, then we
define $\int_0^t\zeta(x,s,\o)\,dw_j(s)$ as the limit  in $Z_t^0$ as
$k\to\infty$ of a sequence $\int_0^t\zeta_k(x,s,\o)\,dw_j(s)$, where
the sequence $\{\zeta_k\}$ is such  as in Proposition \ref{propL}.
\end{definition}
\begin{definition} 
\label{defsolution} {\rm We say that equations
(\ref{parab1})-(\ref{parab10}) are satisfied for $u\in Y^1$ if there
exists $(\chi_1,...,\chi_N)\in Y^1\times (X^0)^N$ such that \baaa
&&u(\cdot,t,\o)=u(\cdot,T,\o)+ \int_t^T\big(\A u(\cdot,s,\o)+
\varphi(\cdot,s,\o)\big)\,ds \ \nonumber
\\&&\hphantom{xxx}+ \sum_{i=1}^N
\int_t^TB_i\chi_i(\cdot,s,\o)ds-\sum_{i=1}^N
\int_t^T\chi_i(\cdot,s)\,dw_i(s)
\label{intur} \eaaa for all $r,t$ such that $0\le r<t\le T$, and
this equality is satisfied as an equality in $Z_T^{-1}$.}
\end{definition}
Note that the condition on $\p D$ is satisfied in the  sense that
$u(\cdot,t,\o)\in H^1$ for a.e. \ $t,\o$. Further, $u\in Y^1$, and
the value of  $u(\cdot,t,\o)$ is uniquely defined in $Z_T^0$ given
$t$, by the definitions of the corresponding spaces. The integrals
with $dw_i$ in (\ref{intur}) are defined as elements of $Z_T^0$. The
integral with $ds$ in (\ref{intur}) is defined as an element of
$Z_T^{-1}$. In fact, Definition \ref{defsolution} requires for
(\ref{parab1}) that this integral must be equal  to an element of
$Z_T^{0}$ in the sense of equality in $Z_T^{-1}$.
\par
In the case where  $\dcoer =0$, Condition \ref{cond3.1.A} is too weak to
ensure solvability of problem (\ref{parab1})-(\ref{parab2}) in
$Y^1$. Therefore, we will need a relaxed version of solution that
does not require $Y^1$-type regularity of $u$.
\subsubsection*{Solution in the representation sense}
For simplicity, we assume in this section that $\varphi\equiv 0$.
\par
Without a loss of generality, we assume that there exist functions
${\ww\b_i: Q\times \O \to \R^n}$, $i=1,\ldots, M$, such that $$
2b(x,t,\o)=\sum_{i=1}^N\b_i(x,t,\o)\,\b_i(x,t,\o)^\top
+\sum_{j=1}^M\,\ww\b_j(x,t,\o)\,\ww\b_j(x,t,\o)^\top, $$ and $\ww
\b_i$ has the similar properties as $\b_i$. (Note that, by Condition
\ref{cond3.1.A}, $2b\ge \sum_{i=1}^N\b_i\b_i^\top$).
\par
 Let
$\ww w(t)=(\ww w_1(t),\ldots, \ww w_M(t))$ be a new Wiener process
independent on $w(t)$.  Let $s\in[0,T)$
be given. Consider the following Ito equation
\begin{eqnarray}
\label{yxs} &&dy(t) =
f(y(t),t)\,dt+\sum_{i=1}^N\b_i(y(t),t)\,dw_i(t) +\sum_{j=1}^M\ww
\b_j(y(t),t)\,d \ww w_j(t),
\nonumber\\ [-6pt] &&y(s)=x.
\end{eqnarray}
\par
Let  $y(t)=y^{x,s}(t)$ be the solution of (\ref{yxs}), and let
$\tau^{x,s}\defi\inf\{t\ge s:\ y^{x,s}(t)\notin D\}$.
\par
\par
To proceed further, we have to impose more conditions.
\par
Let $r_1,r_2\in\R$ be such that $r_1<r_2$.
\par
For the case where  $n=1$, set $ O\defi\{x\in\R: r_1< x< r_2\}$. For
the case where  $n>1$, we assume that $r_1>0$ and  $
O\defi\{x\in\R^n: r_1< |x|< r_2\}$, i.e.,  it is a spherical
layer.
\par
We assume that the following condition is satisfied.
\begin{condition}\label{condO} At least one of the following conditions is satisfied.
\begin{enumerate}
\item The functions $f(x,t)$, $\b_i(x,t)$, $\ww\b_i(x,t)$ are non-random; or
\item
There exists a bijection $\phi:D\to O$ such that the
process $\w y^{x,s}(t)\defi \phi(y^{x,s}(t))$ is such that, for any
$(x,s)\in Q$, there exist bounded random processes $\w
f^{x,s}:[0,+\infty)\times\O\to\R$ and
$h^{x,s}(t)=(h_1^{x,s}(t),....h_{N+M}^{x,s}(t)):[s,+\infty)\times\O\to\R^{N+M}$
that are progressively measurable with respect to $\F_t$ such that
\baaa
\d_h\defi \inf_{x,s}\essinf_{t,\o}|h^{x,s}(t,\o)|^2>0
\eaaa and
\baaa\label{(1r)}dr^{x,s}(t)=\w f^{x,s}(t)dt+\sum_{k=1}^N\h^{x,s}_k(t) dw_k(t)+\sum_{k=N+1}^{N+M}\h^{x,s}_k(t)d\ww w_{k-N}(t), \eaaa where  $r^{x,s}(t)\defi |\w y^{x,s}(t)|$, $t\ge s$.  \end{enumerate}
\end{condition}
\par
Condition \ref{condO}(ii) covers the following two cases:
\begin{itemize}
\item  $n=1$, and $D$ is a connected interval; or
\item  $n>1$, and $D\defi D_1\backslash D_0$, where
$D_i\subset \R^n$ are domains with $C^2$-smooth boundaries $\p D_i$,
$i=0,1$, such that $D_0\subset D_1$ and $\p D_0\cap\p
D_1=\emptyset$. In other words,  there is a lacuna $D_0$ in the
domain $D$.
\end{itemize}
Note that, in  both cases, there exists a bijection $\phi:D\to O
$ such that $\phi$ is continuously twice differentiable inside $D$,
and the derivatives are uniformly bounded. The verification of the
conditions required is straightforward. \par
 For $t\ge s$,
set
\baa
\g^{x,s}(t) \defi\exp\left(-\int_s^t
\lambda(y^{x,s}(t),t)\,dt\right).\label{gamma}
\eaa
\begin{definition} 
\label{defsolU} We say that differential equation (\ref{parab1})
with $\varphi=0$  is satisfied for $u\in U$ in the representation
sense if, for any $(x,s)\in Q$, \baa \hbox{the process}\quad
\g^{x,s}(t\land \tau)u(y^{x,s}(t\land \tau),t\land
\tau)\quad\hbox{is a martingale}. \label{repM} \eaa\end{definition}
\begin{remark}{\rm Definition \ref{defsolU} allows to consider
solutions of differential equation  (\ref{parab1}) without any
requirements on their differentiability.
  }\end{remark}
  \par
 A justification for this definition is the following. First,
 assume that $\xi\in\V$ is given and $\G=0$. In this case,
$u\in U$ is uniquely defined by (\ref{repM}) and (\ref{parab2}),
since it follows from these equations that \baa
u(x,s)=\E\{\g^{x,s}(T\land \tau)\xi(y^{x,s}(T\land \tau),T\land
\tau)|\F_s\}. \label{rep}\eaa In addition, (\ref{parab10}) holds for
any $u\in U$. Second, property (\ref{repM}) holds for the
traditional solution from Definition \ref{defsolution}; this can be
seen from the following.
\begin{theorem}\label{ThJust}  Assume that Condition \ref{cond3.1.A} holds with some
$\dcoer >0$. Let $\xi\in\V\cup Z_T^0$,     and let equations
(\ref{parab1})-(\ref{parab10}) with $\varphi=0$ be satisfied for
 $u\in Y^1$  in the sense of Definition \ref{defsolution} with $\chi_i\in\X^0$.
 Then  (\ref{repM}) holds for this $u$. In other words,  equation (\ref{parab1}) is satisfied for this $u$
 in the representation sense.
\end{theorem}
 \begin{remark}{\rm
Alternatively,   solution of boundary problem
(\ref{parab1})-(\ref{parab2}) could be defined directly by
(\ref{rep}) without requiring (\ref{repM}). \index{However, this does not
allow to describe the set of solutions of differential equation (\ref{parab1})
without specifying the boundary conditions.} We prefer
Definition \ref{defsolU} since it allows to consider  differential equation
(\ref{parab1}) separately from the boundary conditions.
}\end{remark}
\section{Backward SPDEs with the standard terminal condition} In
this section, we assume that $\G=0$. In this case, condition
(\ref{parab2}) can be rewritten as \baa u(\cdot,T)=\xi.
\label{parabC} \eaa
\begin{lemma}
\label{lemma1}  Assume that Condition \ref{cond3.1.A} holds with
$\dcoer >0$. Let $k=0$ or $k=1$. Then problem
(\ref{parab1})-(\ref{parab10}), (\ref{parabC}) has a unique
solution $(u,\chi_1,...,\chi_N)$ in the class
$Y^{k+1}\times(X^1)^N$  for any $\varphi\in X^{k-1}$ and $\xi\in
Z_T^k$, and \be \label{4.2} \| u
\|_{Y^{k+1}}+\sum_{i=1}^N\|\chi_i\|_{X^k}\le C
(\|\varphi\|_{X^{k-1}}+\|\xi\|_{Z^k_T}), \ee where $C>0$ does not
depend on $\xi$.
\end{lemma}\par
For $k=1$, this result is well known; see, e.g., Dokuchaev (1992) or
Theorem 4.2 from Dokuchaev (2010). For $k=2$,  Lemma \ref{lemma1} is
a reformulation of Theorem 3.1 from Du and Tang (2012), or Theorem
3.4 from Dokuchaev (2011), or Theorem 4.3 from Dokuchaev (2012) (the
preprint of this paper was web-published in 2006).  Note that in
Dokuchaev (2011, 2012)  some strengthened version of Condition
\ref{cond3.1.A} was required (Condition 3.5 in Dokuchaev (2011) or
equivalent Condition 4.1 in Dokuchaev (2012)).  The result  in Du
and Tang (2012) was obtained by different methods without this
restriction, i.e., under Condition \ref{cond3.1.A} only.
\par
The following theorem covers the most difficult degenerate case, i.e., where
 Condition \ref{cond3.1.A} holds with $\dcoer =0$ only.
\begin{theorem}
 \label{ThRep}  For any $\xi\in\V$,
there is a unique $u\in U$ such that equations
(\ref{parab1})-(\ref{parab10}) are satisfied in the representation
sense, and that (\ref{parabC}) is satisfied as an equality in
$Z^0_T$. In addition, \baaa
\label{EstRep} \|u\|_{U}\le C_\lambda\|\xi\|_{\V}, \eaaa where
$C_\lambda=\exp(T\sup_{x,t,\o}\max(0,\lambda(x,t,\o)))$.
\end{theorem}
\section{Backward SPDEs with a non-local boundary condition}
In this section, we assume that the following conditions are satisfied.
\begin{condition}\label{condL}
$\lambda(x,t,\o)\le 0$ a.e.
\end{condition}
\begin{condition}\label{condG} The mapping $\G: U_{PC}\to \V$ is linear and continuous and such that at least
one of the following condition holds:
\begin{enumerate}
\item
 $\|\G u\|_{\V}\le \|u\|_{U}$ for any
$u\in U$, and that there exists $\t<T$ such that $\G u=\G
(\Ind_{\{t\le \t\}}u)$.
\item There exists $c\in (0,1)$ such that
 $\|\G u\|_{\V}\le c\|u\|_{U}$ for any
$u\in U$.
\end{enumerate}
\end{condition}
\begin{example} {\rm
Condition \ref{condG}(i) is satisfied for the following operators:
\begin{enumerate}
\item $\G u=\kappa u(\cdot,0)$, $\kappa\in[-1,1]$;
\item
$ (\G u)(x,\o)=\kappa u(x,t_1,\o),\quad t_1\in[0,T);$
\item
$(\G u)(x,\o)=\zeta(\o) u(x,t_1,\o),\quad t_1\in[0,T),\qquad
\zeta\in L_{\infty}(\O,\P,\F_T,\P),\quad |\zeta(\o)|\le 1\quad
\hbox{a.s.} $;
\item
$ (\G u)(x,\o)=\a_1 u(x,t_1,\o)+\a_2 u(x,t_2,\o),\quad
t_1,t_2\in[0,T),\quad |\a_1|+|\a_2|\le 1$;
\item \baaa (\G
u)(x,\o)=\int_0^{\t}k(t)u(x,t,\o)dt,\quad \t\in[0,T),\qquad
k(\cdot)\in L_{\infty}(0,\t),\qquad  \int_0^\t |k(t)|dt\le 1; \eaaa
  \item
\baaa (\G u)(x,\o)=\int_0^{\t}dt\int_Dk(t,y,x,\o)u(y,t,\o)dy, \eaaa
where $\t\in[0,T)$, $k(\cdot):[0,\t]\times D\times D\times\O$ is a
bounded measurable function from $L^{\infty}
(\O,\F_T,\P,L_{\infty}([0,\t]\times D\times D))$
 such that \baaa \esssup_{(x,\o)\in D\times\O}\int_0^\t dt\int_D |k(t,x,y,\o)|dy\le 1. \eaaa
\end{enumerate}
Convex combinations of the operators from this list are also
covered.}
\end{example}
\begin{example} {\rm
Condition \ref{condG}(ii) is satisfied for the following operators:
\begin{enumerate}
\item \baaa (\G
u)(x,\o)=\int_0^{T}k(t)u(x,t,\o)dt,\qquad k(\cdot)\in
L_{\infty}(0,\t),\qquad  \int_0^T |k(t)|dt< 1; \eaaa
  \item
\baaa (\G u)(x,\o)=\int_0^{T}dt\int_Dk(t,y,x,\o)u(y,t,\o)dy, \eaaa
where $\t\in[0,T)$, $k(\cdot):[0,T]\times D\times D\times\O$ is a
bounded measurable function from $L^{\infty}
(\O,\F_T,\P,L_{\infty}([0,\t]\times D\times D))$
 such that \baaa \esssup_{(x,\o)\in D\times\O}\int_0^T dt\int_D |k(t,x,y,\o)|dy< 1. \eaaa
\end{enumerate}
Convex combinations of the operators from this list are also
covered.}\end{example}\par Additional examples of admissible $\G$ can be found is Section
\ref{SecF} below.
\begin{theorem}\label{Th6} For any $\xi\in \V$, there exists a unique $u\in U$ such
that equations (\ref{parab1})-(\ref{parab10})     with $\varphi\equiv 0$ are satisfied in the
representation sense, and \baa \|u\|_{U}\le C\|\xi\|_{\V},
\label{mainest}\eaa where $C>0$ does not depend on $\xi$.
\end{theorem}
\section{Application: portfolio selection problem}\label{SecF}
Theorem \ref{Th6} can be applied to portfolio selection problem for
 continuous time diffusion market model,
 where the market dynamic is described by stochastic differential equations.
Examples of these models can be found in, e.g., \index{\cite{K}} Ma
and Yong (1997) and Karatzas and   Shreve (1998).
\par
We consider the following stripped to the bone  model of a securities market
consisting of a risk free bond or bank account with the price $B(t)
$, ${t\ge 0}$, and
 a risky stock with price $S(t)$, ${t\ge 0}$. The prices of the stocks evolve
 as \be \label{S} dS(t)=S(t)\left(a(t)dt+\s(t) dw(t)\right), \quad
t>0, \ee where $w(t)$ is a Wiener process, $a(t)$ is an
appreciation rate, $\s(t)$ is a volatility coefficient.
The initial price $S(0)>0$ is a  given deterministic constant. The
price of the bond evolves as \baaa \label{Bond} B(t)=e^{rt}B(0), \eaaa
where $B(0)$ is a  given constant, $r\ge 0$ is a short rate. For
simplicity, we assume that $r=0$   and $B(t)\equiv B(0)$.
\par
 We  assume that $w(\cdot)$ is a standard Wiener process on  a
given standard probability space $(\Omega,\F,\P)$, where $\Omega$ is
a set of elementary events, $\F$ is a complete $\s$-algebra of
events, and $\P$ is a probability measure.
\par
 Let $\F_t$
be the filtration generated by $w(t)$.  In particular, this means that $\F_t$ is independent from
$\{w(t_2)-w(t_1)\}_{t_2\ge t_1\ge t}$, and $\F_0$ is trivial, i.e.,
it is the $\P$-augmentation of the set $\{\emptyset,\Omega\}$.
\par
We assume that the processes  $a(t)$, $\s(t)$, and $\s(t)^{-1}$ are bounded and $\F_t$-adapted and continuous. In
particular, this means that the process the process $a(t)$ can be random.
\index{We assume that the processes  $a(t)$, $\s(t)$, and $\s(t)^{-1}$ are bounded, the process  $(S(t),a(t))$ is $\F_t$-adapted, and  the process $\s(t)$
is continuous and deterministic. In
particular, this means that the process the process $a(t)$ can be random and that the process
$(S(t),a(t))$ is currently
observable.}\subsubsection*{Strategies for bond-stock-options market}
The rules for the operations of the agents on the market define the
class of admissible strategies  where the optimization problems have
to be solved.

Let $X(0)>0$ be the initial wealth at time $t=0$ and let $X(t)$ be
the wealth at time $t>0$.
\par
We assume that the wealth $X(t)$ at time $t\in[0,T]$ is
\begin{equation}
\label{X} X(t)=\b(t)B(t)+\g(t)S(t).
\end{equation}
Here $\b(t)$ is the quantity of the bond portfolio, $\g(t)$ is the
quantity of the stock  portfolio, $t\ge 0$. The pair $(\b(\cdot),
\g(\cdot))$ describes the state of the bond-stocks securities
portfolio at time $t$. Each of  these pairs is  called a strategy.
\par  A pair $(\b(\cdot),\g(\cdot))$  is said to be an admissible
strategy if the processes $\b(t)$ and $\g(t)$ are progressively
measurable with respect to the filtration $\F_t$.

 In particular, the agents are not supposed to know the
future (i.e., the strategies have to be adapted to the flow of
current market information).

In addition, we require that  $$ \E\int_0^{T}\left(\b(t)^2B(t)^2+
S(t)^2\g(t)^2\right)dt<+\infty.$$

A pair $(\b(\cdot),\g(\cdot))$  is said to be an admissible
self-financing strategy, if
 $$
dX(t)=\b(t)dB(t)+\g(t)dS(t).$$
Since $B(t)\equiv B(0)$, this means that
 \baaa
 X(t)=X(0)+\int_0^t\g(s)dS(s),
 \eaaa
 and the process
$\g(t)$ alone defines the strategy.

Let $\P_*$ be an equivalent probability  measure such that $S(t)$ is
a martingale under $\P_*$. By the assumptions on $(a,\s)$, this measure
exists and is unique. Under this measure, $X(t)$ is a martingale as
well.
\subsection*{A special goal achieving problem}
In portfolio theory, a typical problem is creation of a strategy
such that the wealth replicates a given contingent claim. The
structure of these claims can be quite complicated; in particular,
these claims may represent payoffs for the derivatives to be hedged. It will be demonstrated below
that Theorem \ref{Th6} can be applied to replication of certain
exotic contingent claims depending on the past portfolio value.
\par
 Without a loss of generality, we
assume that $\P$ is a  martingale  probability measure, i.e., $S(t)$
is a martingale and $dS(t)=\s(t)S(t)dw(t)$.

Let us consider the following example.
In portfolio theory, a typical problem is creation of a strategy
such that the wealth replicates a given contingent claim. The
structure of these claims can be quite complicated; in particular,
these claims may represent payoffs for the derivatives to be hedged. It will be demonstrated below
that Theorem \ref{Th6} can be applied to replication of certain
exotic contingent claims depending on the past portfolio value.

Let $\t\in(0,T]$,  $s_\LS\in(0,S(0))$, $s_\US\in(S(0),+\infty)$, $W_\LS>0$, $W_\US>0$ be given.
Further, let $k_i(t,\o)$ be given random processes such that $k_i(t)$ are
$\F_T$-measurable for any $t$ and \baaa \int_0^\t(|k_1(t)|+|k_2(t)|)dt\le
1\quad\hbox{a.s.} \eaaa In addition, we require that if $\t=T$ then
there exists $\oo c\in[0,1)$ such that \baaa \int_0^T(|k_1(t)|+|k_2(t)|)dt\le \oo
c\quad\hbox{a.s.} \eaaa

Let $\zeta=\zeta(x,\o)\in\oo\V$ be
given.

\par Let $\tau=\inf\{t: S(t)\notin[s_\LS,s_\US]\}$.
\par
We consider the following goal-achieving problem: find an initial wealth $X(0)$
and a self-financing  portfolio strategy such that the corresponding wealth $X(t)$
is bounded and such that
 \baa
 &&X(t)=W_\LS\quad\hbox{if}\quad t\ge \tau\quad\hbox{and}\quad S(\tau)=s_\LS,\label{XL}\\
&&X(t)=W_\US\quad\hbox{if}\quad t\ge \tau\quad\hbox{and}\quad S(\tau)=s_\US,\label{XU}\\
&&X(T)=\int_0^{\t} k_1(s)
X(s)ds+\E\int_0^{\t}k_2(s)X(s)ds+\zeta(S(T))\quad\hbox{if}\quad
\tau>T.\label{XT} \eaa
We impose below some additional restrictions on the choice of $\zeta$
to ensure that $X(T)$ and $X(\tau)$ are similarly defined in the where$\tau=T$.

This toy example still has an economic
meaning. In (\ref{XL})-(\ref{XT}), $\tau$ is the liquidation time
when portfolio is converted into cash. Conditions
(\ref{XL})-(\ref{XU}) may describe preferences for the case of  the
extreme event $\tau<T$; this can be considered as an extreme event
if $s_\LS$ is sufficiently small and $s_\US$ is sufficiently large.
Condition (\ref{XT}) can be illustrated by the following examples.
Let $k_1(t)$ and let  $k_2(t)$ be some $\F_t$-adapted
processes
representing the dividend payoff rate. The manager wishes that the wealth after deduction of the dividends
meets certain target represented by $\zeta$ that  takes into account the dividend payments and the expected deductions caused by the dividend  payments.

In other examples, $k_i(t)$ may represent the annual proportion rate for
the proportional hedge management reward.   Condition (\ref{XT})
with large positive $\zeta$ may represent a goal of a hedge fund
manager who wishes to demonstrate a strong growth expressed via a strong dominance of the resulting wealth
over the total amount of the paid management fees  (for $k_2(t)\equiv 0$), or over the expected
 total amount of the paid management  fees for $k_1(t)\equiv 0$).

Our setting covers also a modification of the
previous examples where  $k_i(t)$ are $\F_T$-measurable, i.e., they are
selected at the terminal time $T$ (say, during annual shareholders meeting),
and the corresponding rates will be applied backward, with respect to the curve of the past wealth.

Let $\oo U={\rm B}([0,T];\oo\V)\cap X^0\cap \CC^0$.

Let   $\ell(x)=c_1x+c_0$, where  $c_1,c_0\in\R$ be such that $\ell(s_\LS)=W_\LS$ and $\ell(s_\US)=W_\US$.
Let $U_{\LU}=\{u:\ u(x,t,\o)=\oo u(x,t,\o)+\ell(x),\quad \oo u\in U\}$. Sometimes, we will consider $\ell$ as an element of the spaces $\oo\V$
or $\oo U$.

For $\eta_\LS,\eta_\US\in L_\infty(\O,\F_T,\P)$, let $\V(\eta_\LS,\eta_\US)=\{v\in\oo\V:\ v(s_\LS)=\eta_\LS,\quad v(s_\US)=\eta_\LS\ \hbox{a.s}\}$.
Here the set $\oo\V$ is defined as above with $n=1$ and $D=(s_\LS,s_\US)$. Let $\V_{\LU}=\V(W_\LS,W_\US)$.

Let $\kappa_i=\int_0^\t
k_i(t)dt$.

 Consider mapping $\G:\oo U\to\oo V$ defined as
\baaa &&
(\G u)(x,\o)\\&&=\int_0^\t
k_1(t)\E\{ u\left(S^{x,t}(T\land\tau),T\land\tau\right)|\F_t\}dt+\E\int_0^\t
k_2(t)\E\{ u\left(S^{x,t}(T\land\tau),T\land\tau\right)|\F_t\}  dt. \eaaa

By the definitions,
\baaa
&&\G u\in \V((\kappa_1+\kappa_2)W_\LS,(\kappa_1+\kappa_2)W_\US),\\
 &&u(\cdot,T)-\G u\in \V( (1-(\kappa_1+\kappa_2)W_\LS, (1-(\kappa_1+\kappa_2)W_\US)\quad\hbox{if}\quad   u\in U_{\LU}.
 \eaaa
Moreover, if $\zeta\in \V(\eta_\LS,\eta_\US)$ for some $\eta_\LS,\eta_\US\in L_\infty(\O,\F_T,\P)$, and $u\in U_{\LU}$, then  \baaa
 &&\zeta+\G u\in \V( \eta_\LS+(\kappa_1+\kappa_2)W_\LS, \eta_\US+(\kappa_1+\kappa_2)W_\US).
 \eaaa
\par
Starting from now, we impose additional restrictions on $\zeta$: we assume that
\baa
\zeta\in \V( (1-\kappa_1-\kappa_2)W_\LS, (1-\kappa_1-\kappa_2)W_\US).\label{zeta}\eaa This condition ensures
that $X(T)$ and $X(\tau)$ are similarly defined in the where$\tau=T$. In addition, this condition ensures   that
 \baaa
 &&\zeta+\G u\in \V_{\LU}\quad\hbox{if}\quad u\in U_{\LU}.
 \eaaa

Let $\xi=\zeta+\G\ell-\ell$, i.e., \baaa
&&\xi(x,\o)=\zeta(x)+\int_0^\t
k_1(t)\ell(x)dt+\E\int_0^\t
k_2(t)\ell(x)dt-\ell(x)\\&&=\zeta(x,\o)+\int_0^\t
k_1(t)\E\{ \ell\left(S^{x,t}(T\land\tau)\right)|\F_t\}dt+\E\int_0^\t
k_2(t)\E\{ u\left(S^{x,t}(T\land\tau)\right)|\F_t\}  dt-\ell(x). \eaaa
By the assumption (\ref{zeta}) on $\zeta$, we have that $\xi\in\V$.
\par
Let us consider
the following problem  \baa
&&d_tu(x,t)+\Bigl(\frac{1}{2}\s(t)^2
x^2\frac{\p^2 u}{\p x^2}(x,t)+\s(t)x\frac{d\chi }{dx}(x,t)\Bigr)=\chi(x,t) dw(t),\quad t<T,\hphantom{xxx}\label{f1}\\
&& u(s_\LS,t)=u(s_\US,t)=0,\label{f2}\\&& u(x,T)-(\G
u)(x)=\xi(x).\label{f3} \eaa  Here $x\in (s_\LS,s_\US)$, This is a special case of problem
(\ref{parab1})-(\ref{parab2})  with $n=N=1$,
$D=(s_\LS,s_\US)$,
 \baaa &&\A
v=\frac{1}{2}\s(t)^2 x^2\frac{\p^2 u}{\p x^2}(x), \qquad B_1v=
\s(t)x\frac{dv}{dx}.\eaaa  The assumptions of Theorem \ref{Th6}
are satisfied for this problem. By this theorem, there exists a
unique solution $u\in U$ of
problem   (\ref{f1})-(\ref{f3})  such that  equation
(\ref{f1}) is satisfied in the representation sense.  In particular,
this means that $u(x,0)=\E\{ u(S(T\land\tau),T\land\tau)|S(0)=x\}$ and
$\E\{u(S(T\land\tau),T\land\tau)|\F_t\}=u(S(t\land \tau),t\land\tau).$
Let  \baaa H(x,t,\o)=u(x,t,\o)+\ell(x). \eaaa
It follows from the definitions that $H(x,T)-(\G
H)(x)=\zeta(x)$.
 \begin{theorem}\label{ThF}
The investment problem  (\ref{XL})-(\ref{XT}) has a solution with the
 wealth  $X(t)=H(S(t\land\tau),t\land\tau)$, $t\in[0,T]$.
\end{theorem}
\index{\begin{remark}{\rm  The statement of Theorem \ref{ThF} holds also for the limit case where
$\t\to 0$, i.e., where (\ref{XT}) is replaced by the constrains
\baaa X(T)=X(0)+\zeta(S(T)), \quad \tau>T.\label{XT0}\eaaa \index{In this
case, the goal achieving problem (\ref{XL})-(\ref{XT}) can be
reduced to replication  of a European barrier option.}  }
\end{remark}}
\begin{remark} {\rm Our approach does not allow to extend Theorem \ref{ThF}
on the case where  $\t=T$ and $\oo c=1$  at the same time. }
\end{remark}
\begin{remark} By the linearity of $\ell$, it follows that
$\E\{H(S(T\land\tau),T\land\tau)|\F_t\}=H(S(t\land \tau),t\land\tau).$ This allows to
extend the definition of the representation solution on the problem for $v$ with non-homogenuous boundary condition on $\p D$  \baaa
&&d_tH(x,t)+\Bigl(\frac{1}{2}\s(t)^2
x^2\frac{\p^2 v}{\p x^2}(x,t)+\s(t)x\frac{d\chi }{dx}(x,t)\Bigr)=\chi(x,t) dw(t),\quad t<T,\hphantom{xxx}\label{v1}\\
&& H(s_\LS,t)=W_\LS, \quad H(s_\US,t)=W_\US,\label{v2}\\&& H(x,T)-(\G
H)(x)=\zeta(x).\label{v3} \eaaa
\end{remark}
\section{Proofs}
 For the brevity, we will use notations $\P_s(\cdot)\defi
\P(\cdot|\F_s)$ and  $\E_s(\cdot)\defi \E(\cdot|\F_s)$.\par We need
the following auxiliary lemma.
\begin{lemma}\label{propnu}
\begin{enumerate}
\item
For any $\vartheta>0$, there exists  $\nu=\nu(\vartheta)\in(0,1)$
that depends only on $D,\A,B_j$ and such that $\P_s(\tau^{x,s}>
s+\vartheta)\le \nu$ a.s. for all $s\ge 0$, and for any $x\in D$.
\item
For any $\vartheta>0$, $\P_s(\tau^{x,s}>s+\vartheta)\to 0$  in
$L_\infty(\O,\F,\P)$ as $\dist(x,\p D)=\inf_{y\in \p D}|x-y|\to 0$.
\item
$\P_s(\tau^{x,s}<\oo s)\to 0$ in $L_\infty(\O,\F,\P)$ as $\oo s-s\to
0$ for any $x\in D$, $s<\oo s<T$.
\end{enumerate}
\end{lemma}
\par {\it Proof of Lemma
\ref{propnu}.} Assume that Condition \ref{condO}(i) is satisfied. In this case, the process $y^{x,s}(t)$ is a Markov diffusion process.
Statement (i) is a reformulation of Lemma 2.1 from Dokuchaev (2004). Statements (ii)-(iii) are well known for  Markov diffusion processes. Hence
Lemma
\ref{propnu} holds if Condition \ref{condO}(i) is satisfied. \par
 Assume that Condition \ref{condO}(ii) is satisfies. We will follow the approach from
 Dokuchaev (2004), p. 296.
\par
Let $D_r\defi (r_1,r_2)$. For $(x,s)\in D\times [0,T)$, we  have that $\tau^{x,s}=\inf\{t\ge
s: \ r^{x,s}(t)\notin D_r\}$ and \be\label{dd1}
\P_s(\tau^{x,s}>s+\vartheta )=\P_s(r^{x,s}(t)\in D_r\ \forall
t\in[s,s+\vartheta ]). \ee
\par
\def\MM{M^{x,s}}\def\BM{B^{x,s}}
Let \baaa \MM(t)\defi
\sum_{k=1}^N\int_s^th_k^{x,s}(r)dw_k(r)+\sum_{k=N+1}^{N+M}\int_s^th_k^{x,s}(r)d\ww
w_k(r),\quad t\ge s. \eaaa  By Condition \ref{condO}(ii), we have
\be\label{dd3}h^{x,s}(t)^\top h^{x,s}(t)= |h^{x,s}(t)|^2\ge
\d_h>0.\ee  Clearly, $\MM(t)$ is a martingale  conditionally
given $\F_s$ vanishing at $t=s$ with quadratic variation process
$$[\MM]_t= \int_s^t|h(y^{x,s}(r),r)|^2dr,\qquad t\ge s.
$$
\par Let $\t^{\mu}(t)\defi \inf\{r\ge s:\ [\MM]_r>t-s\}$.
 Note that $\t^{\mu}(s)=s$, and
the function $\t^{\mu}(t)$ is  strictly increasing in $t>s$ given
$(x,s)$.  By Dambis--Dubins--Schwarz Theorem (see, e.g., Revuz and
Yor (1999)), the process $\BM (t)\defi \MM  (\t^{\mu}(t))$ is a
Brownian motion  conditionally given $\F_s$ vanishing at $t=s$, i.e.,
$\BM (s)=0$, and $\MM(t)=\BM (s+[\MM]_t)$.
\par
Let us prove statement (i).  Let $\w D_r\defi (r_1+K_1,r_2+K_2)$,
where \baaa
K_1\defi -r_2-\vartheta \sup_{x,s,t,\o}|\w f^{x,s}(t,\o)|,\quad
K_2\defi -r_1+\vartheta \sup_{x,s,t,\o}|\w f^{x,s}(t,\o)|.\eaaa It is
easy to see that \be\label{dd2} \P_s(r^{x,s}(t)\in D_r\ \forall
t\in[s,s+\vartheta ])\le \P_s(\MM(t)\in \w D_r\ \forall
t\in[s,s+\vartheta ]). \ee
Clearly, \baa \label{dd4}
\P_s(\MM(t)\in \w D_r \hphantom{x}\forall
t\in[s,s+\vartheta ])\nonumber=\P_s(\BM (s+[\MM]_t)\in \w D_r
\hphantom{x}\forall
t\in[s,s+\vartheta ])\nonumber\\
\le \P_s(\BM (q)\in \w D_r\hphantom{x}\forall
q\in[s,s+[\MM]_{s+\vartheta }]).\eaa By (\ref{dd3}),
$[\MM]_{s+\vartheta }\ge c_2 \vartheta $ a.s. for all $x,s$.
Hence \be\label{dd5}\P_s(\BM (q)\in \w D_r\ \hphantom{x}\forall
r\in[s,s+[\MM]_{s+\vartheta }])\le\P_s(\BM (q)\in \w D_r\
\hphantom{x}\forall q\in[s,s+\d_h   \vartheta ]).\hphantom{xxxxx} \ee By
(\ref{dd1})--(\ref{dd2}) and (\ref{dd4})--(\ref{dd5}), it follows
that
$$
\sup_{x,s}\P_s(\tau^{x,s}>s+\vartheta )\le \nu, $$
where \baaa
\nu\defi \sup_{x,s}
\P_s(\BM (q)\in \w D_r\ \hphantom{x}\forall q\in[s,s+\d_h  \vartheta
]).
\eaaa
We have that $\nu\in(0,1)$ and it depends on $D,\A,B_j$ only.
 This completes the
proof of Lemma \ref{propnu}(i).\par
Let us prove statement (ii).
Clearly, $\dist(x,\p D)\to 0$ if and only if $|\phi(x)|\to r_1$ or
$|\phi(x)|\to r_2$.
\par
Assume that $|\phi(x)|\to r_2$.
 Let $\ww D_r\defi (r_1-KT,r_2)$, where
$K\defi \sup_{x,s,t,\o}|\w f^{x,s}(t,\o)|$. We have that \baa
\P_s(r^{x,s}(t)\in D_r\ \forall t\in[s,s+\vartheta ])\hphantom{xxxxxxxxxx}\nonumber\\\le
\P_s(|\phi(x)|+\MM(t)-K(t-s)\in\ww D_r\ \forall t\in[s,s+\vartheta
]).\label{dd2z} \eaa Clearly,
\baa
\P_s(|\phi(x)|+\BM (s+[\MM]_t)-K(t-s)\in\ww D_r\hphantom{x}\forall
t\in[s,s+\vartheta ])\hphantom{xxxxx}\nonumber\\
\le \P_s(\BM (q)-K\t^\mu(q)\in [r_1-r_2-KT,r_2-|\phi(x)|]\
\hphantom{x}\forall q\in[s,s+[\MM]_{s+\vartheta }]).\label{dd4z}
\eaa By (\ref{dd3}), $[\MM]_{s+\vartheta }\ge \d_h  \vartheta $
a.s. for all $x,s$. Hence \baa \P_s(|\phi(x)|+\BM (q)-K\t^\mu(q)\in
D_r\ \hphantom{x}\forall q\in[s,s+[\MM]_{s+\vartheta }])\hphantom{xxx}\nonumber\\
 \le\P_s(\BM (q)-K\t^\mu(q)\in [r_1-r_2-KT,r_2-|\phi(x)|] \hphantom{x}\forall q\in[s,s+\d_h   \vartheta ]).\label{dd5z} \eaa By (\ref{dd1}),
(\ref{dd2z}),(\ref{dd4z}), and (\ref{dd5z}), and by the properties of a Brownian motion, it follows that
$$
\P_s(\tau^{x,s}>s+\vartheta )\to 0\quad\hbox{in}\quad
L_\infty(\O,\F,\P)\quad \hbox{as}\quad |\phi(x)|\to r_2.
$$
The case where $|\phi(x)|\to r_1$ can be considered similarly.
 This completes the
proof of Lemma \ref{propnu}(ii).
\par
Let us prove statement (iii). Let $\oo D_r(t)\defi (r_1+K(t-s),r_2-K(t-s))$, where $K\defi
\sup_{x,s,t,\o}|\w f^{x,s}(t,\o)|$, and where $t$ is close enough to
$s$ such that $r_1+K(t-s)<r_2-K(t-s)$. We have that \baa
\P_s(r^{x,s}(t)\in D_r\ \forall t\in[s,s+\vartheta ])\hphantom{xxx}\nonumber \\ \ge
\P_s(|\phi(x)|+\MM(t)\in\oo D_r(t)\ \forall t\in[s,s+\vartheta
]).\label{dd2zz} \eaa  Clearly,  \baa
\P_s(|\phi(x)|+\BM (s+[\MM]_t)\in\ww D_r(t)\
\hphantom{x}\forall
t\in[s,s+\vartheta ])\hphantom{xxx}\nonumber\\\ge \P_s(\BM (q)-K\t^\mu(q)\in \oo D_r(\t^\mu(q))
\hphantom{x}\forall q\in[s,s+[\MM]_{s+\vartheta
}]).\label{dd4zz} \eaa
By (\ref{dd3}),  $[\MM]_{s+\vartheta }\in [\d_h  \vartheta,C_h\vartheta]$ a.s. for all $x,s$ for some $C_h>\d_h$. Hence
  \baa \P_s(|\phi(x)|+\BM (q)-K\t^\mu(q)\in
D_r\ \hphantom{x}\forall q\in[s,s+[\MM]_{s+\vartheta }])\hphantom{xxxx}\nonumber\\
 \ge\P_s(\BM (q)-K\t^\mu(q)\in [r_1-r_2-KT,r_2-|\phi(x)|] \hphantom{x}\forall q\in[s,s+C_h   \vartheta ]).\label{dd5zz} \eaa   By (\ref{dd1}),
(\ref{dd2zz}), (\ref{dd4zz}), and (\ref{dd5zz}), it follows that  \baaa
\P_s(\tau^{x,s}>s+\vartheta )\to 1\quad\hbox{in}\quad
L_\infty(\O,\F,\P)\quad \hbox{as}\quad \vartheta\to 0+.\hphantom{xxxx}
\eaaa  This completes the
proof of Lemma \ref{propnu}.
 $\Box$
\par {\it Proof of Theorem \ref{ThJust}}. For the case where
$u\in\X_c^2$ and $\chi_j\in\X_c^1$, this theorem follows immediately
from the proof of Lemma 4.1 from Dokuchaev (2011).
\begin{remark}\label{remDu} { The results in  Dokuchaev (2011) were stated under some more restrictive condition than
Condition \ref{cond3.1.A} with $\dcoer >0$ (Condition 3.5 in the cited
paper). Thanks  to Theorem 3.1. from Du and Tang (2012), this
additional condition can be lifted, i.e., all results from Dokuchaev
(2011) are still valid if  Condition 3.5 from this paper is replaced
by Condition \ref{cond3.1.A}.}
\end{remark}
\par
Let us consider the general case. We introduce operators
   $$ \A^* v\defi\sum_{i,j=1}^n\frac{\p^2 }{\p x_i \p x_j}
\left(b_{ij}(x,t)v(x)\right)-\sum_{i=1}^n\frac{\p}{\p x_i }\left(
f_i(x,t)v(x)\right)+\lambda(x,t)v(x)$$ and \baaa B_i^*v\defi
 -\sum_{k=1}^n \frac{\p }{\p x_k }\,\big(\beta_{ik}(x,t,\o)\,v(x)),\qquad i=1,\ldots ,N. \label{AB*}\eaaa Here $b_{ij}$, $x_i$,
$\b_{ik}$ are the components of $b$,  $\b_i$, and $x$.
\par
Let $\rho\in Z_{s}^0$, and let $p=p(x,t,\o)$ be the solution of the
problem \baaa &&d_tp=\A^* p\, dt +
\sum_{i=1}^NB^*_ip\,dw_i(t), \quad t\ge s,\nonumber\\
&&p|_{t=s}=\rho,\quad\quad p(x,t,\o)|_{x\in \p D}=0.\label{p}\eaaa
By Theorem 3.4.8 from Rozovskii (1990), this boundary value problem
has a unique solution $p\in Y^1(s,T)$. Introduce an operator $\M_s:
Z_s^0\to Y^1(s,T)$  such that $p=\M_s\rho$, where $p\in Y^1(s,T)$ is
the solution of this boundary value problem.
\par
 Let $\rho\in Z_s^0$ be such that $\rho\ge 0$ a.e. and $\int_{D}\rho(x)dx=1$ a.s.
 Let $a\in L_2(\O,\F,\P;\R^n)$ be
  such that $a\in D$ a.s. and it
has the conditional probability  density function $\rho$ given
$\F_s$. We assume that $a$ is independent from  $(w(t_1)-w(t_0),\ww
w(t_1)-\ww w(t_0)\}$, $s<t_0<t_1$. Let $p=\M_s\rho$, and let
  $y^{a,s}(t)$ be
the solution of  Ito equation (\ref{yxs}) with the initial condition
$y(s)=a$.
\par
To prove the theorem, it suffices to show that \baaa
\g(t\land\tau^{x,s})u(y^{x,s}(t\land\tau^{x,s}),t\land\tau^{x,s})=\E_t\g(T\land\tau^{x,s})
u(y^{x,s}(T\land\tau^{x,s}),T\land\tau^{x,s})\quad \hbox{a.s.}\eaaa for any $t$. For this, it suffices to prove that \baa
&&\E\int_D\rho(x)
\g(t\land\tau^{x,s})u(y^{x,s}(t\land\tau^{x,s}),t\land\tau^{x,s})\nonumber\\&&=\E\int_D\rho(x)\E_t
\g(T\land\tau^{x,s})u(y^{x,s}(T\land\tau^{x,s}),T\land\tau^{x,s})\label{ident}\eaa
 for any $\rho\in
Z_s^0$ such as described above.
\par
By Theorem 6.1 from Dokuchaev (2011) and Remark \ref{remDu}, we have that \baaa
\int_Dp(x,t)u(x,t)dx=\E_t
\g^{a,s}(t\land\tau^{a,s})u(y^{a,s}(t\land\tau^{a,s}),t\land\tau^{a,s})\eaaa
and \baaa \int_Dp(x,T)u(x,T)dx=\E_T\,\g(T\land\tau^{a,s})
u(y^{a,s}(T\land\tau^{a,s}),T\land\tau^{a,s}).\eaaa By the duality
established in Theorem 3.3 from Dokuchaev (2011) and Remark \ref{remDu}, it follows that
\baaa\E\int_Dp(x,t)u(x,t)dx=\E\int_Dp(x,T)u(x,T)dx.\eaaa
This means that $\E(\E_tq(a,s,t))=\E(\E_Tq(a,s,T))$,
where $$q(a,s,t)=\g^{a,s}(t\land\tau^{a,s})u(y^{a,s}(t\land\tau^{a,s}),t\land\tau^{a,s}).$$
Hence \baa
\E(\E_tq(a,s,t))=\E(\E_tq(a,s,T)).\label{eqeq}\eaa
\par
Without a loss of generality, we shall assume that $a$ is a random vector on the probability space
$(\ww\O,\ww\F,\ww\P)$, where $\ww\O=\O\times\O'$, where $\O'=D$,     $\ww\F=\overline{\F_s\otimes {\cal B}_D}$,
where ${\cal B}_D$ is the set of Borel subsets of $D$, and
\baaa
\ww\P(S_1\times S_2)=\int_{S_1}\P(d\o)\P'(\o,S_2),\quad \P'(\o,S_2)=\int_{S_2}\rho(x,\o)dx
\eaaa
for $S_1\in\F$ and $S_2\in {\cal B}_D$. The symbol $\ww\E$ denotes the expectation in $(\ww\O,\ww\F,\ww\P)$. We suppose that
$\ww\o=(\o,\o')$, $\ww\O=\{\oo\o\}$, and $a(\ww w)=\o'$.
\par
We have that
\baaa
\E(\E_t q(a,s,t))=\E\int_{\ww\O}\ww\P(d\o|\F_t)q(\o',s,t,\o)=\E\int_Dd\o'\rho(\o')\int_{\O}\P(d\o|\F_t)q(\o',s,t,\o)\\=\E\int_D\E_t\rho(\o') q(\o',s,t,\o)d\o'=\E\int_D\rho(\o')q(\o',s,t,\o)d\o'\\
=\E\int_D\rho(x) \g^{x,s}(t\land\tau^{x,s})u(y^{x,s}(t\land\tau^{x,s}),t\land\tau^{x,s})dx\eaaa
and
\baaa
\E(\E_tq(a,s,T))=\E\int_{\ww\O}\ww\P(d\o|\F_t)q(\o',s,T,\o)=\E\int_Dd\o'\rho(\o')\int_{\O}\P(d\o|\F_t)q(\o',s,T,\o)
\\=\E\int_D\rho(\o')\E_t q(\o',s,T,\o)d\o'
=\E\int_D\rho(x)
\g^{x,s}(T\land\tau^{x,s})\E_t u(y^{x,s}(T\land\tau^{x,s}),T\land\tau^{x,s})dx.
\eaaa  Since the choices of $\a$ and $\rho$ are arbitrary, it
follows from (\ref{eqeq})  that (\ref{ident}) holds. This completes
the proof of Theorem \ref{ThJust}. $\Box$
\begin{lemma}\label{lemmaV} Let $\xi\in \V$, and let $u$ be
defined by (\ref{rep}). Then $u(\cdot,s)\in\V$ for any $s$.
\end{lemma}
\def\oox{\oo x}
\def\oos{s}
\par
{\em Proof of Lemma \ref{lemmaV}}.
  By Theorem II.8.1 from Krylov
(1980) applied on the conditional probability space given $\F_s$,
\index{p.143(rus),} we have that \baa
\E_s\sup_{t\in[s,T]}|y^{\oox,s}(t)-y^{x,s}(t)|^2\to 0\quad
\hbox{as}\quad \oox\to x.\label{Kr} \eaa  In addition, we have that
\baa \E_s|\ttau^{\oox,s}-\ttau^{x,s}|\to 0\quad \hbox{as}\quad
\oox\to x \quad\hbox{a.s.}, \label{tau}\eaa where
$\ttau^{x,s}=T\land\tau^{x,s}$. If Condition \ref{condO}(i) is
satisfied, then (\ref{tau}) follows from Theorem 2.3 from Dokuchaev
(2004). If Condition \ref{condO}(ii) is satisfied, then (\ref{tau})
follows from Theorem 2 from Dokuchaev (2008) applied on the
conditional probability space given $\F_s$.

Clearly,  $\u(\cdot,s)\in L_\infty(\O,\F_s,P,B(\oo D,\R))$.
\def\gg{\g^{x,s}(\ttau^{x,s})}
\def\ggg{\g^{\oox,\oos}(\ttau^{\oox,\oos})}
 By the definitions,
$u(x,s)=\E_s\gg\xi(y^{x, s}(\ttau^{x,s}))$. Hence \baaa
&&|u(\oox,\oos)-u(x,s)|=|\E_{\oos}\ggg\xi(y^{\oox,\oos}(\ttau^{\oox,\oos}))-\E_s\gg\xi(y^{x,s}(\ttau^{x,s}))|\\
&&\le  \E_{
\oos}\ggg|\xi(y^{\oox,\oos}(\ttau^{\oox,\oos}))-\xi(y^{x,s}(\ttau^{x,\oos}))|+\E_{
\oos}|\ggg-\gg|\,|\xi(y^{x,s}(\ttau^{x,s})|.\eaaa
 By (\ref{Kr}),(\ref{tau}), and by the Lebesgue Dominated
Convergence Theorem, we obtain that $u(\oox,s)\to u(x,s)$ a.s. Hence $u(\cdot,s)\in\oo\V$.
\par
Further, \baaa u(x,s)=\E_s\gg\xi(y^{x, s}(\ttau^{x,s}))=
\E_{s}\Ind_{\{ \tau^{x,s}>T\}}\gg\xi(y^{x,s}(T)).\eaaa Hence \baaa
|u(x,s)|\le [\E_s\Ind_{\{ \tau^{x,
s}>T\}}^2]^{1/2}[\E_{s}\gg^2\xi(y^{x,s}(T))^2]^{1/2}\\=
(\P_s(\tau^{x,
s}>T)^{1/2}[\E_{s}\gg^2\xi(y^{x,s}(T))^2]^{1/2}
\eaaa By Lemma \ref{propnu}(iii), we have that $\P_s(\tau^{x,
s}>T)\to 0$ a.s. as $\dist(x,\p D)\to 0$. Hence
$u(x,s)\to 0$ as $\dist(x,\p D)\to 0$.
 This
completes the proof of Lemma \ref{lemmaV}. $\Box$
\begin{lemma}\label{lemmaRep}
Let $\xi\in \V$, and let $u$ be defined by (\ref{rep}). Then
(\ref{repM}) holds for this $u$.
\end{lemma}
\par
{\em Proof of Lemma \ref{lemmaRep}}.   Let $\xi\in\V$ be given, and
let $u$ be defined by (\ref{rep}). It suffices to show that
(\ref{repM}) holds. Let $\w B\defi 2 b-\sum_{j=1}^M \b_j\b_j^\top$.
By Condition \ref{cond3.1.A}, the matrix $\w B=\w B(x,t,\o)$ is
non-negatively defined for all $(x,t,\o)$. Let $y_\d^{x,s}(t)$ and
$\g_\d^{x,s}(t)$  be defined similarly to $y^{x,s}(t)$ and
$\g^{x,s}(t)$  such that the corresponding $\ww\b_j$ and $\w B$ are such that\index{bylo: $\w \b$?} are
such that $\sum_{j=1}^M \ww\b_j\ww\b_j^\top\equiv\w B+\d I$, where $I$
is the unit matrix in $\R^{n\times n}$. Let
$\tau_\d^{x,s}=\inf\{t>s:\ y_\d^{x,s}(t)\notin D\}$ and
$\tttau^{x,s}=T\land \tau_\d^{x,s}$.

By Theorem II.8.1 from Krylov (1980)  applied on the conditional
probability space given $\F_s$, \index{p.143(rus),} \baa
\E_s\sup_{t\in[s,T]}|y^{x,s}(t)-y_\d^{x,s}(t)|^2\to 0\quad
\hbox{as}\quad \d\to 0. \label{Kr2}\eaa In addition, \baa
\E_s|\ttau^{x,s}-\tttau^{x,s}|\to 0\quad \hbox{as}\quad \d\to 0.
\label{tau2}\eaa  If Condition \ref{condO}(i) is satisfied, then
(\ref{tau2} follows from Theorem 2.3 from Dokuchaev (2004). If
Condition \ref{condO}(ii) is satisfied, then (\ref{tau2}) follows
from Theorem 2 from Dokuchaev (2008)  applied on the conditional
probability space given $\F_s$.
\par
Let $u_\d$ be defined by (\ref{rep}) with $y^{x,s}(t)$ and
$\tau^{x,s}$ replaced by $y_\d^{x,s}(t)$ and $\tau_\d^{x,s}$
respectively.  By  (\ref{Kr2})-(\ref{tau2}), it follows that, for
any $(x,s)$, there exists a sequence $\d=\d_i\to 0$ such that \baa
\sup_{t\in[s,T]}|y^{x,s}(t)-y_\d^{x,s}(t)|\to 0,\quad
|\ttau^{x,s}-\tttau^{x,s}|\to 0 \quad \hbox{a.s.}\quad
\hbox{as}\quad \d_i\to 0  \label{Kr2as}\eaa and \baa
|y^{x,s}(\ttau^{x,s})-y_\d^{x,s}(\tttau^{x,s})|\to 0,\quad
|\g^{x,s}(\ttau^{x,s})-\g_\d^{x,s}(\tttau^{x,s})|\to 0 \quad
\hbox{a.s.}\quad \hbox{as}\quad \d_i\to 0.  \label{tau2as} \eaa

Let us show that, for all $s$, \baa u_\d(x,s)\to u(x,s) \quad
\hbox{a.s. for all  $x$ as}\quad \d=\d_i\to 0. \label{uulim}\eaa
\par
\def\gg{\g^{x,s}(\ttau^{x,s})}
\def\ggg{\g_\d^{x,s}(\tttau^{x,s})}

Let $(x,s)\in Q$ be given. Let \baaa
z(t)=\g^{x,s}(t\land\tau^{x,s})u(y^{x, s}(t\land\tau^{x,s})),\quad
z_\d(t)=\g_\d^{x,s}(t\land\tau_\d^{x,s})u_\d(y_\d^{x, s}(t\land
\tau_\d^{x,s})).\eaaa
 By the definitions,
\baaa u(x,s)=\E_s\gg\xi(y^{x, s}(\ttau^{x,s}))=\E_sz(T),\quad
u_\d(x,s)=\E_s\ggg\xi(y_\d^{x, s}(\tttau^{x,s}))=\E_sz_\d(T).\eaaa
Hence \baaa |u_\d(x,s)-u(x,s)| =|\E_{s}z_\d(T)-\E_s z(T)| &\le& \E_{
\oos}\ggg|\xi(y_\d^{x,\oos}(\tttau^{x,\oos}))-\xi(y^{x,s}(\ttau^{x,\oos}))|\\&+&
\E_{ \oos}|\ggg-\gg|\,|\xi(y^{x,s}(\ttau^{x,s})|.\eaaa By
(\ref{Kr2as})-(\ref{tau2as}) and by the Lebesgue Dominated
Convergence Theorem,  (\ref{uulim}) holds.
\par
Further, let us estimate the value  \baa \E\int_s^T |z_\d(t)-
z(t)|dt =\ww\psi_1+\ww\psi_2,\label{psipsi} \eaa where \baaa
\ww\psi_1=\E\int_s^{\tau^{x,s}\land \tau_\d^{x,s}\land T} |z_\d(t)-
z(t)| dt,\quad \ww\psi_2=\E\int_{\tau^{x,s}\land \tau_\d^{x,s}\land
T}^{(\tau^{x,s}\lor \tau_\d^{x,s})\land T} |z_\d(t)- z(t)| dt.
\label{eps1e}\eaaa We have that \baaa \ww\psi_1=
\E\int_s^{\tau^{x,s}\land \tau_\d^{x,s}\land
T}|\g^{x,s}_\d(t)u_\d(y_\d(t),t)-\g^{x,s}(t)u(y(t),t)|dt\to
0\quad\hbox{as}\quad \d=\d_i\to 0, \eaaa and\baaa \ww\psi_2\le
(\esssup_{t,\o}|z_\d(t)|+\esssup_{t,\o}|z(t)|)\E|({\tau^{x,s}\lor
\tau_\d^{x,s}})\land T-{\tau^{x,s}\land \tau_\d^{x,s}}\land T|\to
0\quad\hbox{as}\quad \d=\d_i\to 0. \label{eps1ee}\eaaa The last two
limits hold by  (\ref{Kr2})-(\ref{uulim}) and by the
Lebesgue Dominated Convergence Theorem. Hence expectation
(\ref{psipsi}) converges to zero as $\d=\d_i\to 0$. Since the
processes $z_\d(t)$ and $z(t)$ are uniformly bounded, it follows
that there exists a subsequence $\{\d_k'\}$ of the sequence
$\{\d_i\}$ such that \baa \E\int_s^T |z_\d(t)- z(t)|^2dt\to
0\quad\hbox{as}\quad \d=\d_k'\to 0.\label{limz}\eaa In addition, if
follows from (\ref{Kr2as})-(\ref{tau2as}) that $\ggg\to \gg$ a.s.
 as
$\d=\d_k'\to 0$. Similar to  (\ref{limz}), we obtain that \baaa
\E|z_\d(T)-z(T)|^2=\E|\ggg \xi(y_\d^{x,s}(\tttau^{x,s})) -\gg
\xi(y^{x,s}(\ttau^{x,s}))|^2\to 0\quad\hbox{as}\quad \d=\d_k'\to
0.\label{limzT} \eaaa By Theorem \ref{ThJust}, $z_\d(t)=\E_tz_\d(T)$
for any $t$, i.e., this process is a martingale in $t\in[s,T]$.
Therefore, the limit process $z(t)$  is also a martingale.
This completes the proof of Lemma \ref{lemmaRep}. $\Box$
\begin{lemma}\label{lemmaC}
Let $\xi\in \V$, and let $u$ be defined by (\ref{rep}). Then
$u\in\CC^0$.
\end{lemma}
\par
{\em Proof of Lemma \ref{lemmaC}}.
 Let $\e>0$ be given. Let us show  that
 \baaa
\|u(\cdot,s)-u(\cdot,t)\|_{Z_T^0}\to 0\quad \hbox{as}\quad
|s-t|\to 0.\label{epsTo}\eaaa  By the definitions, \baa
\E|\u(x,s)-u(x,t)|^2\le \Psi_1(x,s,t)+\Psi_2(x,s,t),
\label{uPP}\eaa where\baaa &&\Psi_1(x,s,t)=
\E|\u(x,s)-u(y^{x,s}(t\land \tau^{x,s}),t\land
\tau^{x,s})|^2,\\&&\Psi_2(x,s,t)=\E|u(y^{x,s}(t\land
\tau^{x,s}),t\land \tau^{x,s})-u(x,t)|^2.
 \eaaa
Let $\eta(x,s,t)=u(y^{x,s}(t\land \tau^{x,s}),t\land \tau^{x,s})$.
By Theorem \ref{ThJust}, (\ref{repM}) holds. By the Martingale
Representation Theorem, it follows that \baaa
\eta(x,s,t)=\eta(x,s,s)+\sum_{k=1}^{N}\int_s^{t}\pi_k(x,s,q)dw_k(q)+\sum_{k=N+1}^{N+M}\int_s^{t}\pi_k(x,s,q)d\ww
w_{k-N}(q) \eaaa for some functions $\pi_k(x,s,t):Q\times
[s,T]\times\O\to\R$ that are $\F_t$-adapted and such that \baaa
\sup_{x,s}\E\sum_{k=1}^{N+M}\int_s^{T}\pi_k(x,s,q)^2dt\le
2\|\xi\|_{\V}. \eaaa By (\ref{repM}), it follows that \baaa
\Psi_1(x,s,t)= \E|\E_s\eta(x,s,t)-\eta(x,s,t)|^2\le
\E\sum_{k=1}^{N+M}\int_s^{t}\pi_k(x,s,q)^2dq\to 0\label{Psi1a}\eaaa
as $|s-t|\to 0$ for all $x$. By the Lebesgue Dominated Convergence
Theorem, \baa \|\Psi_1(\cdot,s,t)\|_{Z_T^0}\to 0\quad \hbox{as}\quad
s-t\to 0.\label{Psi1} \eaa
\par
Let us estimate $\Psi_2$.
 Clearly,
\baaa \E_s|y^{z,s}(t\land\tau^{x,s})-x|^2\to 0\quad \hbox{as}\quad
|s-t|\to 0\quad \hbox{for all}\quad x.\label{Kr3} \eaaa By Lemma
\ref{propnu}(iii), \baaa \P_s(t\land\tau^{x,s}=t)\to 1\quad
\hbox{as}\quad |s-t|\to 0\quad \hbox{for all}\quad x.\label{Kr0}
\eaaa Hence $\Psi_2(x,s,t)\to 0$ as $s-t\to 0$ for all $x$. By the
Lebesgue Dominated Convergence Theorem,
$\|\Psi_2(\cdot,s,t)\|_{Z_T^0}\to 0$ as $s-t\to 0$. The proof of
Lemma \ref{lemmaC} follows from this limit and from
(\ref{uPP})--(\ref{Psi1}). $\Box$
\par
 {\em Proof of Theorem \ref{ThRep}} follows immediately from
Lemma \ref{lemmaV}, Lemma \ref{lemmaRep}, and Lemma \ref{lemmaC}. Estimate (\ref{EstRep})
follows from (\ref{rep}).
$\Box$\vspace{0.5cm}
\par
To proceed with the proof of Theorem \ref{Th6}, we need to introduce first  some additional definitions.
Let $s\in (0,T]$, $\varphi\in X^{-1}$ and $\Phi\in Z^0_s$. Consider
the problem \be \label{4.1}
\begin{array}{ll}
d_tu+\left( \A u+ \varphi\right)dt +
\sum_{i=1}^NB_i\chi_i(t)dt=\sum_{i=1}^N\chi_i(t)dw_i(t), \quad t\le s,\\
u(x,t,\o)|_{x\in \p D}, \\
 u(x,s,\o)=\Phi(x,\o).
\end{array}
 \ee
 \par
Assume that $\dcoer =0$ in Condition \ref{cond3.1.A}. Introduce
operators $\L_T :\V\to U$, such that $u=\L_T \Phi,$ where
 $u$ is the solution of  problem (\ref{4.1}) in the representation sense. By Theorem \ref{ThRep}, these linear operators
are continuous.
\par
Introduce   operators $\Q:\V\to \V$ a
such that $\Q \Phi=\G\L_T\Phi$, i.e.,  $\Q\Phi= \G u$,
    where
$u\in U$ is the solution in the representation sense of  problem
(\ref{4.1}) with $\Phi\in \V$. Since the operator $\G:U\to \V$ is
continuous, the operators $\Q:\V\to \V$ is linear and
continuous. In particular, $\|\Q\|\le \|\G\|\|\L_T\|$, where
$\|\Q\|$, $\|\G\|$, and $\|\L_T\|$, are the norms of the operators
$\Q: \V\to \V$, $\G: U\to \V$, and $\L_T: \V\to U$, respectively.
\begin{lemma}\label{lemmaQ} If the operator $(I-\Q)^{-1}:\V\to \V$ is  continuous
then problem (\ref{4.1})   has a unique solution $u\in U$ in the
representation sense for any $\xi \in \V$. For this solution, \baa
\label{Q} u=\L_T (I-\Q)^{-1}\xi\eaa
 and \baaa
\label{Qes} \| u \|_{U}\le C \|\xi\|_{\V}, \eaaa where $C>0$ does
not depend on $\xi$.
\end{lemma}
\par
{\it Proof of Lemma \ref{lemmaQ}}.    Clearly,
 $u\in U$ is the solution of   problem
(\ref{parab1})-(\ref{parab2}) if and only if \baa
&&u=\L_Tu(\cdot,T), \label{Q1}\\
&&u(\cdot,T)-\G u=\xi. \label{Q2}\eaa
 Since $ \G u=\Q
u(\cdot,T)$, equation (\ref{Q2})  can be rewritten as \baa
 u(\cdot,T)-\Q u(\cdot,T)
=\xi. \label{Q3}\eaa  By the continuity of $(I-\Q)^{-1}$, equation
(\ref{Q3}) can be rewritten as
$$ u(\cdot,T)=(I-\Q)^{-1}\xi. $$  Therefore, equations (\ref{Q1})-(\ref{Q2}) imply that \baaa \label{4.3}
u=L_T\varphi+\L_T u(\cdot,T)=\L_T (I-\Q)^{-1}\xi. \eaaa Further, let
us show that if (\ref{Q}) holds then equations (\ref{Q1})-(\ref{Q2})
hold. Let $u$ be defined by (\ref{Q}). Since $u=\L_T u(\cdot,T)$, it
follows that $u(\cdot,T)=(I-\Q)^{-1}\xi$. Hence \baaa u(\cdot,T)-\Q
u(\cdot,T)=\xi,\eaaa i.e., $u(\cdot,T)-\G \L_T u(\cdot,T)=\xi.$
Hence \baaa u(\cdot,T)-\G \L_T u(\cdot,T)=\xi.\eaaa This means that
(\ref{Q1})-(\ref{Q2}) hold. Then the proof of Lemma \ref{lemmaQ}
follows. $\Box$
\par
We are now ready to prove  Theorem \ref{Th6}.
\par
 {\em Proof of Theorem \ref{Th6}}.  Let $\|\Q\|_{\V,\V}$ be
the norm of the operator $\Q=\G\L_T:\V\to \V$.
 By (\ref{rep}), we have  for $u=\L_T\Phi$ that
  \baaa
\sup_{s\in[0,T]}\|u(\cdot,s)\|_{\V}\le
\|\Phi\|_{\V}.\label{estRs}\eaaa Hence \baa  \|\L_T\Phi\|_{\W}\le
\|\Phi\|_{\V}.\label{estR}\eaa By the assumptions on $\G$, it
follows that $\|\G u\|_{\V}\le c\|u\|_{U}$, where $c<1$ if Condition \ref{condG}(ii) is satisfied.
It follows  that if  Condition \ref{condG}(ii) holds then
\baa
\|\Q\|_{\V,\V}\le c<1. \label{Q<1}\eaa
\par
Let us assume that Condition \ref{condG}(i) is satisfied.
\par \index{Assume first that $c_\lambda=0$.that $\w\lambda$ is replaced
by $\ww\lambda=\w\lambda+|c_\lambda|$.}
Let $u=\L_T\Phi$, $s\in[0,T]$. Let $y(t)=y^{x,s}(t)$ be the solution
of  Ito equation (\ref{yxs}) with the initial condition $y(s)=x$.
By (\ref{rep}), it follows that \baaa
&&\|u(\cdot,s)\|_{\V}=\esssup_{x,\o}\E_s\g^{x,s}(T)\Phi(y^{x,s}(T))\Ind_{\{
\tau^{x,s}\ge T\}}\\  &&\le \esssup_{x,\o} \left[\E_s\Ind_{\{
\tau^{x,s}\ge T\}}^2\right]^{1/2} \esssup_{x,\o}
\left[\E_s\Phi(y^{x,s}(T))^2\right]^{1/2}
\\
 &&\le   \esssup_{x,\o} \left[\E_s\Ind_{\{ \tau^{x,s}\ge
T\}}^2\right]^{1/2} \|\Phi\|_{\V}=\esssup_{x,\o} \P_s(\tau^{x,s}\ge
T) ^{1/2} \|\Phi\|_{\V}. \eaaa If $s<\t$ then $\{\tau^{x,s}\ge
T\}\subseteq \{\tau^{x,s}\ge s+\vartheta\}$, where $\vartheta\defi
T-\t>0$. Hence
 \baaa \|u(\cdot,s)\|_{\V}\le\esssup_{x,\o}
\P_s(\tau^{x,s}\ge s+\vartheta) ^{1/2} \|\Phi\|_{\V},\quad s\le \t.
\eaaa
  It follows that \baaa
\|u(\cdot,s)\|_{\V} \le \nu^{1/2} \|\Phi\|_{\V},\quad s\le \t
\eaaa and  \baaa \|\Ind_{\{s\le \t\}}u\|_{\W} \le \nu^{1/2}
\|\Phi\|_{\V}. \eaaa By  Condition \ref{condG}(ii) on $\G$, it follows that
\baaa \|\Q\Phi \|_{\V} =\|\G u\|_{\V}=\|\G (\Ind_{\{s\le \t\}}u)\|_{\V} \le
\nu^{1/2} \|\Phi\|_{\V},\quad s\le \t, \quad u=\L_T\Phi. \eaaa
 It follows that if  Condition \ref{condG}(i) holds then  \baa
 \|\Q\|_{\V,\V}\le  \nu^{1/2}<1.\label{Q<nu}\eaa

 By   (\ref{Q<1}) and (\ref{Q<nu}), it follows that the operator
    $(I-\Q)^{-1}:\V\to \V$ is bounded. Let \baa
\label{uQ} u=\L_T (I-\Q)^{-1}\xi.\eaa By the assumptions on $\G$ and
by (\ref{estR}), it follows that $\xi+\TT\varphi=\xi\in\V\subset
\V$. Hence $(I-\Q)^{-1}\xi\in\V$. By the properties of $\L_T$ , it
follows that $u\in U$.
 Similar to the proof of Lemma \ref{lemmaQ}, it can be shown that $u$ is a
 solution  of problem
(\ref{parab1})-(\ref{parab2}) in the representation sense. Estimate
(\ref{mainest}) follows from the continuity of the corresponding
operators in (\ref{uQ}). Then the proof of Theorem \ref{Th6}
follows. $\Box$
\par
{\em Proof of Theorem \ref{ThF}}. Assume that the terminal
wealth is defined as $X(T)=H(S(T\land\tau),T\land\tau)$. By the
Martingale Representation Theorem, there exists $\psi(t)$ such that
\baaa H(S(\ttau),\ttau)= \E
H(S(\ttau),\ttau)+\int_0^T\psi(t)dw(t)=\E
H(S(\ttau),\ttau)+\int_0^T\g(t)dS(t), \eaaa  where
$\ttau=T\land\tau$ and $\g(t)=\psi(t)S(t)^{-1}\s(t)^{-1}$.
Therefore, $X(t)\defi \E\{H(S(\ttau),\ttau)|\F_t\}$ is the wealth
for the self-financing strategy such that $X(0)=\E
H(S(\ttau),\ttau)$.

By the linearity of $\ell$ and martingale property of $S(t)$, we
have that $$\E\{\ell(S(\ttau))|\F_t\}=\ell(S(t\land \tau)).$$ Since
$u$ is the solution of  problem (\ref{parab1})-(\ref{parab2}) in the
representation sense, we have that
$$\E\{u(S(\ttau),\ttau)|\F_t\}=u(S(t\land \tau),t\land\tau).$$ Hence
$X(t)=H(S(t\land \tau),t\land \tau)$.
\par
By the choice of $u$ and $\ell$, conditions  (\ref{XL})-(\ref{XU})
are satisfied for this $X(t)$.
\par
Let us show that condition (\ref{XT}) is satisfied for this $X(t)$.
If $\tau>T$ then, by the definitions,
\baaa &&X(T)-\int_0^\t
k_1(t)X(t)dt-\E \int_0^\t k_2(t)X(t)dt\\&&=H(S(T),T) -\int_0^\t
k_1(t)H(S(t\land \tau),t\land \tau)dt-\E
\int_0^\t k_2(t)H(S(t\land \tau),t\land \tau)dt\\&&=
H(S(T),T) -\int_0^\t k_1(t)
\E\{H(S(T\land \tau),T\land \tau)|\F_t\}dt\\&&\hphantom{xxxxxxxxxxxxx}-\E\int_0^\t
k_2(t)
\E\{H(S(T\land \tau),T\land \tau)|\F_t\}dt\\&&= H(S(T\land \tau),T\land \tau)-(\G
  u)(S(T\land \tau))=
  \xi(S(T\land \tau))=\xi(S(T)).
\eaaa
 This completes the proof of Theorem \ref{ThF}. $\Box$
\subsection*{Acknowledgment} This work  was
supported by ARC grant of Australia DP120100928 to the
author.
\section*{References} $\hphantom{XX}$
Al\'os, E., Le\'on, J.A., Nualart, D. (1999).
 Stochastic heat equation with random coefficients
 {\it
Probability Theory and Related Fields} {\bf 115} (1), 41--94.
\par
Bender, C. and  Dokuchaev, N. (2014). A
first-order BSPDE for swing option pricing.   {\em
Mathematical Finance}, in press.
\par
Bally, V., Gyongy, I., Pardoux, E. (1994). White noise driven
parabolic SPDEs with measurable drift. {\it Journal of Functional
Analysis} {\bf 120}, 484--510.
\par Caraballo,T., P.E. Kloeden, P.E.,
Schmalfuss, B. (2004). Exponentially stable stationary solutions for
stochastic evolution equations and their perturbation, Appl. Math.
Optim. {\bf  50}, 183--207.
\par
 Chojnowska-Michalik, A. (1987). On processes of Ornstein-Uhlenbeck type in
Hilbert space, Stochastics {\bf 21}, 251--286.

\par
Chojnowska-Michalik, A. (1990). Periodic distributions for linear
equations with general additive noise, Bull. Pol. Acad. Sci. Math.
{\bf 38} (1�12) 23--33.

\par
Chojnowska-Michalik, A., and Goldys, B. (1995). {Existence,
uniqueness and invariant measures for stochastic semilinear
equations in Hilbert spaces},  {\it Probability Theory and Related
Fields},  {\bf 102}, No. 3, 331--356.
\par
Da Prato, G., and Tubaro, L. (1996). { Fully nonlinear stochastic
partial differential equations}, {\it SIAM Journal on Mathematical
Analysis} {\bf 27}, No. 1, 40--55.
\par
Dokuchaev, N.G. (1992). { Boundary value problems for functionals
of
 Ito processes,} {\it Theory of Probability and its Applications}
 {\bf 36} (3), 459-476.
 \par Dokuchaev, N.G. (2004).
Estimates for distances between first exit times via parabolic
equations in unbounded cylinders. {\it Probability Theory and
Related Fields}, {\bf 129}, 290 - 314.
\par
Dokuchaev, N.G. (2005).  Parabolic Ito equations and second
fundamental inequality.  {\it Stochastics} {\bf 77} (2005), iss. 4.,
pp. 349-370.
\par
Dokuchaev, N. (2008a). Estimates for first exit  times of
non-Markovian It\^o processes.  {\it Stochastics} {\bf 80},
397--406.
\par
 Dokuchaev N. (2008b) Parabolic Ito equations with mixed in time
conditions.
{\it Stochastic Analysis and Applications} {\bf 26}, Iss. 3, 562--576. 
\par
Dokuchaev, N. (2010). Duality and semi-group property for backward
parabolic Ito equations. {\em Random Operators and Stochastic
Equations. } {\bf 18}, 51-72.
\par
Dokuchaev, N. (2011). Representation of functionals  of Ito
processes in bounded domains. {\em Stochastics} {\bf 83}, No. 1,
45--66.
\par
 Dokuchaev, N. (2012).
Backward parabolic Ito equations and second fundamental inequality.
{\em Random Operators and Stochastic Equations} {\bf 20}, iss.1,
69-102. Preprint: web published in 2006,
http://lanl.arxiv.org/abs/math/0606595.
\par
Dokuchaev, N. (2013).
On forward and backward SPDEs with non-local boundary conditions.
Accepted to {\em Discrete and Continuous Dynamical Systems --
Series A (DCDS-A)}.
Available at http://xxx.lanl.gov/abs/1307.8223.
\par
Du K., and Tang, S. (2012). Strong solution of backward stochastic
partial differential equations in $C^2$ domains. {\em  Probability
Theory and Related Fields)} {\bf 154}, 255--285.
\par
Duan J.,  Lu K., Schmalfuss B. (2003). Invariant manifolds for
stochastic partial differential equations. {\em Ann. Probab.} {\bf
31}
 2109--2135.
\par
Feng C., Zhao H. (2012). Random periodic solutions of SPDEs via
integral equations and Wiener-Sobolev  compact embedding. {\em
Journal of Functional Analysis} {\bf 262}, 4377--4422.
\par
Gy\"ongy, I. (1998). Existence and uniqueness results for semilinear
stochastic partial differential equations. {\it Stochastic Processes
and their Applications} {\bf 73} (2), 271-299.
\par
Hamza, K.,  Klebaner, F.C. (2005)
On Solutions of First Order Stochastic Partial Differential Equations.
Working paper: http://arxiv.org/abs/math/0510495.
 \par
  Hu, Y., Ma, J., Yong, J. (2002). On semi-linear degenerate backward stochastic partial differential equations.
  Probab. Theory Related Fields. Vol. 123 (3), pp. 381--411.
\par
Karatzas, I., and Shreve, S.E. (1998).   {\em Methods of
Mathematical Finance}, Springer-Verlag, New York.
\par
Kl\"unger, M. (2001). Periodicity and Sharkovsky�s theorem for
random dynamical systems, {\em Stochastic and Dynamics} {\bf 1},
iss.3, 299--338.
Krylov, N.V. (1980).  {\em Controlled diffusion processes}.
Shpringer, New York.
\par Krylov, N. V. (1999). An
analytic approach to SPDEs. Stochastic partial differential
equations: six perspectives, 185--242, Mathematical Surveys and
Monographs, {\bf 64}, AMS., Providence, RI, pp.185-242.
\par
Ladyzhenskaia, O.A. (1985). {\it The Boundary Value Problems of
Mathematical Physics}. New York: Springer-Verlag.
\par
Liu, Y., Zhao, H.Z (2009). Representation of pathwise stationary
solutions of stochastic Burgers equations, {\em Stochactics and
Dynamics} {\bf  9} (4), 613--634.
\par
 Ma, J., Yong, J. (1997). Adapted solution of a class
 of degenerate backward stochastic partial differential equations, with applications.
{\em Stochastic Processes and Their Applications.} Vol. 70, pp.
59-84. \par
 Ma, J., Yong, J. (1999). On linear, degenerate backward
stochastic partial differential equations.
 Probability Theory and Related Fields. Vol. 113 (2), pp. 135--170.
\par Maslowski, B. (1995). { Stability of semilinear equations with
boundary and pointwise noise}, {\it Annali della Scuola Normale
Superiore di Pisa - Classe di Scienze} (4), {\bf 22}, No. 1, 55--93.
\par
Mattingly, J. (1999). Ergodicity of 2D Navier�Stokes equations with random forcing and large viscosity. {\em Comm. Math.
Phys.} 206 (2),  273�288.
\par
Mohammed, S.-E.A., Zhang T.,  Zhao H.Z. (2008). The stable manifold theorem for semilinear stochastic evolution equations
and stochastic partial differential equations. {\em Mem. Amer. Math. Soc.} 196 (917)  1�105.
\par
Pardoux, E. (1993). Stochastic partial differential equations, a
review, Bull. Sci. Math. {\bf 117}, no. 1, 29-47.
\par
 Revuz, D., and Yor, M. (1999). {\it Continuous Martingales
and Brownian Motion}. Springer-Verlag: New York.
\par
Rodkina, A.,  Dokuchaev, N., and Appleby, J. (2014)
On limit  periodicity of discrete time stochastic processes.  {\em
Stochastic and Dynamics}, in press.
\par
Rozovskii, B.L. (1990). {\it Stochastic Evolution Systems; Linear
Theory and Applications to Non-Linear Filtering.} Kluwer Academic
Publishers: Dordrecht-Boston-London.
\par
Sinai, Ya. (1996). Burgers system driven by a periodic stochastic
flows, in: Ito's Stochastic Calculus and Probability Theory,
Springer, Tokyo, 1996, pp. 347--353.
\par
Walsh, J.B. (1986). An introduction to stochastic partial
differential equations, Lecture Notes in Mathematics {\bf 1180},
Springer Verlag.
\par
Yong, J., and Zhou, X.Y. (1999). { Stochastic controls: Hamiltonian
systems and HJB equations}. New York: Springer-Verlag.
\par
 Zhou, X.Y. (1992). { A duality analysis on stochastic partial
differential equations}, {\it Journal of Functional Analysis} {\bf
103}, No. 2, 275--293.
\end{document}